\documentclass[12pt, twoside]{article}
\usepackage{amsmath,amsthm,amssymb}
\usepackage{times}
\usepackage{enumerate}
\usepackage{graphics}
\usepackage{appendix}
\usepackage{subfigure}
\usepackage{tikz}

\def\titlerunning#1{\gdef\titrun{#1}}
\makeatletter
\def\author#1{\gdef\autrun{\def\and{\unskip, }#1}\gdef\@author{#1}}
\def\address#1{{\def\and{\\\hspace*{18pt}}\renewcommand{\thefootnote}{}%
\footnote {#1}}%
\markboth{\autrun}{\titrun}}
\makeatother
\def\email#1{e-mail: #1}
\def\subjclass#1{{\renewcommand{\thefootnote}{}%
\footnote{\emph{Mathematics Subject Classification (2010):} #1}}}
\def\keywords#1{\par\medskip
\noindent\textbf{Keywords.} #1}

\def\TT{{\Bbb T}}
\def\CC{{\Bbb C}}
\def\NN{{\Bbb N}}
\def\ZZ{{\Bbb Z}}
\def\SS{{\Bbb S}}
\def\RR{{\Bbb R}}

\def\L{{\Lambda}}
\def\t{{\theta}}

\def\l{{\lambda}}

\def\pa{{\partial}}

\def\e{{\epsilon}}
\def\beq{\begin{equation}}
\def\eeq{\end{equation}}

\def\qedbox{$\rlap{$\sqcap$}\sqcup$}
\def\skipaline{\removelastskip\vskip12pt plus 1pt minus 1pt}

\def\Proof{\removelastskip\skipaline
\noindent \it Proof. \rm}

\newtheorem{Theorem}{Theorem}
\newtheorem{Lemma}{Lemma}[section]

\newtheorem{Proposition}{Proposition}[section]
\newtheorem{Corollary}{Corollary}[section]
\newtheorem{Remark}{Remark}[section]

\begin{document}


\baselineskip=17pt


\titlerunning{Positive Lyapunov Exponent are not Open}

\title{Quasi-Periodic Schr\"odinger Cocycles with Positive Lyapunov Exponent are not Open in the Smooth Topology}

\author{Yiqian Wang
\and
Jiangong You}

\date{}

\maketitle

\address{Y. Wang:
 Department of Mathematics,
Nanjing University, Nanjing 210093, China; \email{yiqianw@nju.edu.cn}
\and
J. You: Department of Mathematics,
Nanjing University, Nanjing 210093, China; \email{jyou@nju.edu.cn}}



\begin{abstract}
 One knows that the set of  quasi-periodic Schr\"odinger cocycles with positive Lyapunov exponent is open and dense in the analytic topology. In this paper, we  construct
 cocycles  with positive  Lyapunov exponent which can be approximated  by ones with zero Lyapunov exponent in the space of ${\cal C}^ l$ ($1 \le l\le \infty$) smooth quasi-periodic cocycles. As a consequence, the set of  quasi-periodic Schr\"odinger cocycles with positive Lyapunov exponent is not ${\cal C}^ l$ open.  

\end{abstract}

\section{ Introduction and Results}
 Let $X$ be a ${\cal C}^r$ compact manifold,   $T: X\to X$ be  ergodic with a normalized invariant measure $\mu$ and $A(x)$ be
 a $SL(2, \mathbb{R})$-valued function on $X$.   The dynamical system: $(x, w)\to (T(x), A(x)w)$  in $X\times \mathbb{R}^2$ is called a $SL(2, \mathbb{R})$ cocycle (or cocycle for simplicity)
  over the base dynamics $(X, T)$. We will simply denoted it as $(T,  A)$.
If the base system is a rotation on torus, i.e., $X=\TT^m=\RR^m\backslash \ZZ^m$, $T= T_\omega: x\to
x+ \omega$  with rational independent $\omega$, we call $(T_\omega, A)$ a
quasi-periodic cocycle, which is  simply
denoted by $(\omega, A)$.
If furthermore $A(x)=S_{v}(x)$ is of the form
$
S_{v}(x)=\left(
\begin{array}{ll}
v(x) & -1\\
\ \ \ 1& 0\end{array} \right)
$
with $v(x+1)=v(x)$, we call $(\omega, S_{v}(x))$ a quasi-periodic Schr\"odinger
cocycle. \\

 For any $n\in \NN$ and $x\in X$, we denote
$$
A^n(x)=A(T^{n-1}x)\cdots A(Tx)A(x)
$$
and
$$
A^{-n}(x)=A^{-1}(T^{-n}x)\cdots A^{-1}(T^{-1}x)=(A^n(T^{-n}x))^{-1}.
$$
If the base dynamics $(X, T, \mu)$ is fixed, the (maximum) Lyapunov exponent of $(T, A)$ is defined as
$$
L(A)=\lim_{n\rightarrow \infty}\frac{1}{n}\int \log\|A^n(x)\|d\mu:=\lim_{n\rightarrow \infty} \int L_n(A(x))d\mu
\in [0, \infty).
$$
$L(A)$ measures the average growth rate of $\|A^n(x)\|$.

The regularity  and positivity  of the Lyapunov exponent (LE) are the central subjects in dynamical systems. One is also interested in the problem whether or not cocycles with positive LE  are open and dense.  The problems turn out to be very subtle, which depend on the base dynamics $(X, T)$ and the smoothness of the matrix $A$.

Firstly, classical Furstenberg theory \cite{Furstenberg} showed that for certain special linear cocycles
over Bernoulli shifts, the largest
LE is positive under very general conditions. Furstenberg and Kifer \cite{FK} and
Hennion \cite{Hennion} proved the continuity of the largest LE of i.i.d random matrices under a condition of almost
irreducibility.  Kotani \cite{Ko2} showed
that the LE of  Schr\"odinger cocycles $S_{E-v}$ is positive for almost every energy $E$ if the potential $v$ is non-deterministic.
Viana\cite{Viana}
 proved that for any $s>0$, the set of $C^s$ linear cocycles over any
hyperbolic ergodic transformation  contains an open and dense subset of cocycles with nonzero
LE; and the LE is continuous for $SL(2,\mathbb{R})$-cocycles over Markov shifts \cite{MV}. For other related results, see \cite{AV}, \cite{BV} and \cite{VY}.

When the base dynamics is uniquely ergodic (e.g., irrational rotation or skew shift on the torus),  the positivity and continuity of the LE seem to be more sensitive to the smoothness of the matrix-valued function $A(x)$.
The LE was proved to be discontinuous at any non-uniform cocycles in the $C^0$ topology by Furman \cite{Furman} (Continuity at uniform hyperbolic cocyces and cocycles with zero LE is trivial).
Motivated by Ma$\tilde{n}\acute{e}$ \cite{Mane1,Mane2}, Bochi
\cite{[Bo1]} further  proved a stronger result that any non-uniformly hyperbolic $SL(2,\mathbb{R})$-cocycle over  a fixed ergodic system on a compact space can be approximated by cocycles with zero LE  in the $C^0$ topology, which
 shows that any non-uniform cocycle can not be an inner point of cocycles with positive LE in the $C^0$ topology. For further related results, we refer to \cite{Bjerk2}, \cite{BF}, \cite{BochiV}, \cite{Hennion},  \cite{JMavi}, \cite{JMavi2}, \cite{Knill}, \cite{Th}.

On the other hand, there are tremendously many positive results in the analytic topology.
Herman \cite{Herman} introduced the subharmonicity method and showed that the LE of $S_{E- 2\lambda \cos x}$ is positive for  $|\lambda|>1$ and all $E$. Herman also proved the positivity of the LE for trigonometric polynomials if the coupling is large enough. The generalization to arbitrary one-frequency non-constant real analytic potentials was obtained by Sorets and Spencer \cite{SoSp}. Same results for \emph{Diophantine} multi-frequency were established by Bourgain and Goldstein \cite{BoGo} and Goldstein and Schlag \cite{GS}. Zhang \cite{ZZhang} gave a different proof of the results in \cite{SoSp} and applied it to a certain class of analytic Szeg\H o cocycle.  For more references, one can see \cite{Bourgain}, \cite{Eliasson}, \cite{Klein}.

 For the continuity of the LE,  Large Deviation Theorems (LDT) is an important tool, which was first established by Bourgain and Goldstein in \cite{BoGo} for real analytic potentials with \emph{Diophantine} frequencies. In \cite{GS}, Goldstein and
Schlag proved that, by some sharp version of LDT and generalized Avalanche Principle(AP), $L(S_{E-v})$ is
H$\rm\ddot{o}$lder continuous in $E$ if $\omega$ is a Diophantine,
$v(x)$ is analytic and  $L(S_{E-v})>0$. Jitomirskaya, Koslover and Schulteis \cite{JKS} proved the continuity of
the LE for a class of analytic one-frequency quasi-periodic $M(2, \CC)$-cocycles with singularities. We will briefly mention more results along this line at the end of this section.
The continuity of the LE
implies that the cocycles with positive LE are open in  analytic topology. Together with the denseness result by Avila \cite{Avila}, one knows that  the set of quasi-periodic cocycles with positive LE is  open and dense in the analytic topology.

We have seen that the behavior of the LE in the $C^0$ topology  is totally different from its behavior in the analytic topology. The smooth case is more subtle.
 Avila\cite{Avila} proved, among other results, that
 the LE is positive for a dense subset of smooth quasi-periodic cocycles. Recently, with Benedicks-Carleson-Young's method\cite{BC,Young}, the authors \cite{WY} constructed  quasi-periodic
cocycles $(T_\omega, A)$ where $T_\omega$ is an irrational rotation $x\to x+
\omega$ on $\SS^1$ with $\omega$ of bounded type and $A\in {\cal C}^l(\SS^1, SL(2,\mathbb{R}))$,
$0\le l\le \infty$, such that the LE is not continuous at $A$ in the ${\cal C}^l$ topology. Such an example in the Schr\"odinger class is also constructed in \cite{WY}. For $\mathcal{C}^2$ cosine-like potentials, Anderson Localization and the positivity of LE has been established by Sinai \cite{sinai} and  Fr\"ohlich-Spencer-Wittwer \cite{frohlich}, also see Bjerkl\"ov \cite{Bjerk3}.
For the model in \cite{sinai}, Wang and Zhang \cite{Wang-Zhang} showed the continuity of the LE, which implies that  non-uniform quasi-periodic cocycles can be   inner points of smooth quasi-periodic cocycles with positive exponents. An interesting  problem is whether or not
quasi-periodic cocycles with positive exponent are open and dense in the smooth topology as in the analytic topology. As we mentioned before, the denseness follows from the result of Avila \cite{Avila}.
 In this paper, we will prove that, different from the analytic case,  the set of smooth quasi-periodic cocycles with positive exponent are not open  in smooth topology.

The LE of quasi-periodic Schr\"odinger cocycles have attracted so much attention  not only because of its importance in dynamical systems, but also due to its close relation with quasi-periodic Schr\"odinger operators. The latter has strong background in physics.
The LE of  Schr\"odingier cocycles coming from  the eigenvalue equations of quasi-periodic Sch\"odinger operators  encodes enormous information on the spectrum.   It is known from Kotani theory that  positive LE implies  singular spectrum, and typically Anderson localization, see \cite{ishii,kotani2,pastur};
while zero Lyapunov spectrum usually implies continuous, typically absolutely continuous spectrum.  The positivity of the LE is also a starting point for many other problems in dynamical systems and spectral theory, such as h\"older continuity of LE, continuity and  topological structure of spectrum set. The recent developed methods, such as  Green's function estimates
and Avalanche Principle, etc.(see \cite{Bourgain0}), depend crucially on  the positivity of the LE.

Another related interesting question is the robustness of Anderson localization. i.e., wether or not the perturbations of a Schr\"odinger operator exhibiting Anderson localization still have  Anderson localization? The answer is affirmative in the analytic category since the LE is continuous and thus the positivity of the LE is kept  under perturbations.\footnote{More precisely, it is true for all almost all frequencies.}
We are interested in the question in smooth case, which is closely related to the problem whether or not the positivity of the LE  is kept under perturbations in the smooth category, equivalently whether or not there exist  smooth non-uniformly hyperbolic Schr\"odinger cocycles which can be accumulated by ones with zero LE in ${\cal C}^l$ topology ($l=1, 2, \cdots, \infty$). If it is the case,  the nature of the spectrum of Schr\"odinger operators might exhibit dramatically changes under small perturbations of the potential in smooth topology.

 The following is the main result of this paper.
\begin{Theorem} \label{T2}
 Consider quasi-periodic Schr\"odinger cocycles over $\mathbb{S}^1 $ with  $\omega$ being a fixed irrational number of
bounded-type.\footnote{ Bounded type means $\frac{p_k}{q_k}$, the best approximation of $\omega$, satisfies $q_{k+1}\le M q_k$ for some $M>0$.} For any $0\le l\le \infty$,
there exists a Schr\"odinger cocycle $S_v$ with  arbitrarily large Lyapunov exponent and a sequence of Schr\"odinger cocycles $S_{v_n}$ with zero Lyapunov exponent such that $v_n(x)\to v(x)$ in the $C^l$ topology.
   As a consequence, the set of quasi-periodic Schr\"odinger cocycles with positive Lyapunov exponent is not $C^l$  open.
 \end{Theorem}

Theorem \ref{T2} can be obtained from Theorem \ref{maintheorem} in the same way as  in \cite{WY} to derive examples in Schr\"odinger cocycles from examples in $SL(2,\mathbb{R})$ cocycles. Thus we only need to prove Theorem \ref{maintheorem}.

\begin{Theorem}\label{maintheorem}
 Consider quasi-periodic $SL(2, \mathbb{R})$ cocycles over $\mathbb{S}^1 $ with  $\omega$ being a fixed irrational number of
bounded-type.  For any $0\le l\le \infty$, there exists a cocycle
$D_l \in {\cal C}^l(\SS^1, SL(2, \mathbb{R}))$ with arbitrarily large Lyapunov exponent and a
sequence of cocycles $C_k\in {\cal C}^l(\SS^1, SL(2, \mathbb{R}))$ with zero Lyapunov exponent such that $C_k\rightarrow D_l$ in the $C^l$ topology. As a consequence, the set of $SL(2, \mathbb{R})$-cocycles with positive Lyapunov exponent is not $C^l$ open.
\end{Theorem}

\begin{Remark} Completely different from the result in Theorem \ref{T2}, Bonatti, G\'omez-Mont and Viana \cite{BGV} proved that there exist H\"older continuous cocycles over Bernoulli shift with positive
LE which can be approximated by continuous cocycles with zero
LE, but not by  H\"older ones, which shows that the base dynamics plays an important role in the regularity problem of the LE.
\end{Remark}
\vskip 0.3cm
\begin{Remark} Avila and Krikorian \cite{AK} showed that  the LE is smooth in the space of smooth monotonic quasi-periodic cocycles. Our result  shows that the monotonicity assumption in \cite{AK} is necessary, and   behavior of the LE in smooth quasi-periodic Schr\"odinger cocycles  homotopic to the identity are completely different from its behavior in the class of monotone cocycles.
\end{Remark}


 The proof of Theorem \ref{maintheorem} is constructive. Recall in \cite{WY}, we have constructed a smooth  cocycles  $D_l$ with positive LE and  a smooth cocycle  $A_1$  in $\frac 1{2k}$-neighborhood of $D_l$  in the  $C^l$ topology for any given $k>0$  such that
 the finite LE of $A_1$, defined by  $L_{n_1}(A_1)=\frac{1}{n_1}\int_{\mathbb{S}^1} \log\|A_1^{n_1}(x)\|dx$, is smaller that $(1-\delta_2)L(D_l)$  for a fixed number $\delta_2>0$.
As a consequence of subadditivity of finite LE,  $L(A_1)< (1-\delta_2) L(A)$. It follows  that the LE is discontinuous at $D_l$. However, the construction in  \cite{WY}  did not tell us how small $L(A_1)$ can be. In this paper we will define a new $A_1$ somehow different from the one in \cite{WY} but satisfies the same property stated as above. Then we further
locally modify $A_1$ such that the modified cocycle, say $A_2$, satisfies $\|A_2-A_1\|_{C^l}<\frac 1{4k}$ and $L_{n_2}(A_2)<(1-\delta_2) L_{n_1}(A_1)$. It follows that $A_2$ is in the $\delta$-neighborhood of $A$ and
$L(A_2)< (1-\delta_2)^2L(A)$. Inductively,  we locally modify $A_i$ such that the modified cocycle, say  $A_{i+1}$, satisfies $\|A_{i+1}-A_i\|_{C^l}<\frac 1{2^ik}$ and $L_{n_{i+1}}(A_{i+1})<(1-\delta_2) L_{n_i}(A_i),$ where
$n_i\to \infty$ will be specified later.  It follows that all $A_{i}$  are in the $\frac 1k$-neighborhood of $D_l$ and
$L(A_{i+1})< (1-\delta_2)^iL(D_l)$. It is easy to see that $A_i$ has a limit, say   $C_k$,  with  $L(C_k)=0$. Moreover, $\|C_k-D_l\|_{C^l}<\frac 1k$. Theorem \ref{maintheorem} is thus proved since $k$ is arbitrary.

 We remark that $D_l$ and $ C_k$ we constructed are  of the forms
$\L R_{\phi(x)}$ and $\L R_{\phi_k(x)}$ where $\L =diag\{\l, -\l\},\l\gg 1$
with $L(D_l)\sim \ln \lambda$ and $L(C_k)=0$. Moreover, $\phi_k(x)$ is an arbitrarily small modification of $\phi(x)$ in an arbitrarily small neighborhood of two special points (called critical points). So a small change makes a big difference! For Schr\"odinger cocycles,  we actually construct, for arbitrarily large but fixed $\l$, smooth $v(x)$ and $\bar v(x)$ which are arbitrarily close to each other and slightly different only at the neighborhood of two critical points such that $L(S_{\l v(x)})$ is very big while $L(S_{\l \bar v(x)})=0$. The result is surprising as we have even not seen any example of smooth Schr\"odinger cocycles of the form  $S_{\l\bar v(x)}$ with
 $\l\gg 1$ such that $L(S_{\l \bar v(x)})=0$.

From our  construction, one can see how and where to modify a cocycle so as to control the LE. This might be useful for other problems.

\noindent
{\it More results on the continuity of the LE in the analytic topology.}
When  the
base dynamics is a shift or skew-shift of a higher
dimensional torus, the log-continuity of the LE was proved in
\cite{BGS} by Bourgain,  Goldstein and Schlag.
Recently,
the result of \cite{JKS} was generalized by Jitomirskaya and Marx \cite{JM1} for all non-trivial singular analytic quasiperiodic cocycles with one-frequency with application to the extended Harper's model \cite{JM2}.

An arithmetic version of large deviations and inductive scheme were developed by Bourgain and Jitomirskaya  in \cite{[BJ]} allowing to obtain joint continuity of the LE for $\mathrm{SL}(2,\mathbb{C})$ cocycles, in frequency and cocycle map, at any irrational frequencies. This result has been crucial in many further important achievements, such as the proof of the Ten Martini problem \cite{AJ}, Avila's global theory of one-frequency cocycles \cite{avila3, avila3.2}. It was extended to multi-frequency case by Bourgain \cite{Bourgain} and to general $\mathrm{M}(2,\mathbb{C})$ case by Jitomirskaya and Marx  \cite{JM2}.
More recently, a completely different method without LDT or AP was developed by Avila, Jitomirskaya and Sadel \cite{avilajitosadel} and was applied to prove the continuity of the LE in $\mathrm{M}(d,\mathbb{C})$, $d\ge 2$.
For further works, see \cite{BoSc}, \cite{duarteklein}, \cite{JM3},  \cite{Klein}, \cite{YZ}.\\

\section{The construction of $D_l$}

We consider the case $m=1$. We say a $SL(2,\mathbb{R})$-matrix $A$ is hyperbolic if $\|A\|> 1.$
A quasi-periodic cocycle $(\omega, A(x))$ of degree $d$ is defined  by a matrix function  $A(x)=R_{\psi(x)}\cdot \L(x)\cdot R_{\phi(x)}$ on $\RR$,  with $\L(x+1)=\L(x)=diag\{\|A\|, \frac1{\|A\|}\} $, $\psi(x+1)=2\pi d+\psi(x), \phi(x+1)=2\pi d+\phi(x)$ where  $R_{\theta}=\left(\begin{array}{ll} \cos \theta\! &\!\! -\sin\! \theta\\
\sin \theta & \cos \theta\end{array}\right)$. It is easy to see that $(\phi(x)+\psi(x-\omega)) $ is uniquely determined by $A(x)$ up to $2\pi \ZZ$ and  $L(A)= L(\L(x)\cdot R_{\phi(x)+\psi(x-\omega)})$ as $A$ is conjugated to   $\L(x)\cdot R_{\phi(x)+\psi(x-\omega)}$. \\

 Let  $\L=diag \{\l, \frac 1\l\}$ with $\l\gg 1$. In this section, we will construct a sequence of smooth cocycles $B_k$ of the form $\L \cdot R_{\xi_k(x)}$, converging in $C^l$ such that
$L(\lim B_k)>0$. Moreover $\xi_k(x)$ will be specially designed so that, in the next section, we can further constructed  cocycles $C_k$ with zero Lyapunov exponent in any small neighborhood of $B_k$. When $\lambda$ is big, we will see that the Lyapunov exponent of $B_k$
 crucially depends on the local behavior, more precisely the degeneracy,  of $\xi_k(x)$ at the critical points
 $\{c: \ \xi_k(c) =\frac\pi 2 \ ({\rm mod}\ \pi)\}$ due to the cancelation.
 The construction in this section is in principle along the line of the construction in \cite{WY}, the difference is in this paper, we use the decomposition of a matrix instead of the most expended and contracted direction of a matrix which makes the proof more transparent.

 Let $\omega$ be a fixed irrational number  and $\frac{p_k}{q_k}$ be its best approximation. Throughout the paper, we assume that $\omega$ is of the bounded type, i.e., $q_{k+1}\le M q_k$;
   $\e>0$ is small. $l$ is a fixed positive integer reflecting the smoothness of cocycles. Let $\l$ and $N$ are large enough  so that \beq \label{N}
  10l\sum_{k=N}^{\infty}
\frac{\log q_{k+1}}{q_k}\le \e, \quad \lambda^{-1}\ll q_N^{-2}. \eeq

We  define the decaying sequence $\{\l_{k}\}$ inductively by $\log\l_{k}=\log\l_{k-1}-\frac{10l\log q_{k}}{q_{k-1}}$ where $\l_N=\l\gg 1$. It is easy to see that
 $\l_k$ converges to $\l_{\infty}$ with $\l_{\infty}>\l^{1-\e}$.
\vskip 0.3cm

 For $k\ge N$, let ${\cal C}_0=\left\{0, \frac12\right\}$, $I_{k,1}=[-\frac{1}{q_{k}^2},\frac{1}{q_{k}^2}]$,
$I_{k,2}=[\frac12-\frac{1}{q_{k}^2},\frac12+\frac{1}{q_{k}^2}]$ and
$I_{k}=I_{k, 1}\bigcup I_{k, 2}$. For $C\ge 1$, we denote by
$\frac{I_{k,1}}{C}=[-\frac{1}{Cq_k^2},
\frac{1}{Cq_k^2}], \frac{I_{k,2}}{C}=[\frac12-\frac{1}{Cq_k^2},
\frac12+\frac{1}{Cq_k^2}],$ and by $\frac{I_k}{C}$ the set
$\frac{I_{k,1}}{C}\cup\frac{I_{k,2}}{C}$.
Denote Lebesgue measure of $I_k$ by $|I_k|$.
For each $x\in I_k$, let $r^+_k(x)$ (respectively $r^-_k(x)$) be the
smallest positive integer such that $T^{r^+_k(x)}({x})\in {I_k}$ (respectively $T^{-r^-_k(x)}({x})\in {I_k}$).   Let $r^{\pm}_k=\min_{x\in I_k} r^{\pm}_k(x)$ and $r_k=\min\{r^+_k, r^-_k\}$. Obviously,
${r_{k}}\ge q_k$. Moreover, from the symmetry between $I_{k,1}$ and $I_{k,2}$, we have $r_k=r_k^+=r_k^-$.

\vskip 0.3cm


  \vskip 0.6 cm




We define $\xi_0$ on $I=I_1\bigcup I_2=\{x: |x|\le \frac1{2q_N^2}\}\bigcup
\{x: |x-\frac12|\le \frac 1{2q_N^2}\}$ by
\begin{equation}\label{phi0}
\xi_{0}(x)=\left\{\begin{array}{ll} \xi_{01}(x),& \ \ |x|\le \frac 1{2q_N^2};\\
 -\xi_{02}(x)\ ({\rm\ or\ } \xi_{02}(x)),& \ \ |x-\frac12|\le
 \frac1{2q_N^2} \end{array} \right.
\end{equation}
where
\begin{equation}\label{phi01}\xi_{01}(x)={\rm sgn}(x)|x|^{l+1},\quad  \xi_{02}(x)={\rm sgn}(x-\frac 12)|x-\frac 12|^{l+1}.
\end{equation}
 $\xi(x)$ is a lift of a $1$-periodic $C^l$ function satisfying
 \begin{equation}\label{phi}
\xi(x)=\left\{\begin{array}{ll} \xi_{01}(x),& \ \ |x|\le \frac 1{2q_N^2};\\
 -\xi_{02}(x)\ ({\rm\ or\ } \pi+\xi_{02}(x)),& \ \ |x-\frac12|\le
 \frac1{2q_N^2}, \end{array} \right.
\end{equation}
and $|\xi(x) ({\rm mod}\ \pi)|> \frac1{2q_N^2}$
 for any $x ({\rm mod}\ 1)\notin I$.
 See Figures 1 and 2 for the picture of $\xi(x)$.

\begin{center}
\begin{tikzpicture}[yscale=1.0]
\draw[thick] (0,0)..controls (0.25,0) and (0.3,0.07)..(0.4,0.2);
\draw[thick](0.4,0.2)..controls (0.5,0.3) and (0.65,0.59)..(0.8,0.6);
\draw[thick] [->] (-1.5,0.0) -- (4,0.0);
\draw [thick][->] (0,-1.5) -- (0,2.5);
\draw[thick] (0.8,0.6)..controls (0.95,0.59) and (1.1,0.3)..(1.2,0.2);
\draw[thick](1.2,0.2)..controls (1.3,0.07) and (1.35,0)..(1.6,0);
\draw[thick] (0+1.6,0)..controls (0.25+1.6,0) and (0.3+1.6,-0.07)..(0.4+1.6,-0.2);
\draw[thick](0.4+1.6,-0.2)..controls (0.5+1.6,-0.3) and (0.65+1.6,-0.59)..(0.8+1.6,-0.6);
\draw[thick] (0.8+1.6,-0.6)..controls (0.95+1.6,-0.59) and (1.1+1.6,-0.3)..(1.2+1.6,-0.2);
\draw[thick](1.2+1.6,-0.2)..controls (1.3+1.6,-0.07) and (1.35+1.6,0)..(1.6+1.6,0);
\draw [dotted, thick] (0.,1) -- (3.5,1);
\node [below] at (-0.31,1.1) {$\pi$};
\draw [dotted, thick] (0.,-1) -- (3.5,-1);
\node [below] at (-0.31,-0.9) {$-\pi$};
\node [below] at (-0.2,0.1) {$0$};
\node [below] at (1.6-0.2,0.1) {$\frac12$};
\node [below] at (3.2,0.1) {$1$};
\node [below] at (4.2,0.2) {$x$};
\node [below] at (0,3.1) {$\xi(x)$};
\node[thick] [below] at (0.5,-1.5) {Figure 1: homotopic to identity};
\end{tikzpicture}\ \ \ \
\begin{tikzpicture}[yscale=1.0]
\draw[thick] (0,0)..controls (0.25,0) and (0.3,0.07)..(0.4,0.1);
\draw[thick](0.6+0.6,0.9)..controls (0.7+0.6,0.93) and (0.75+0.6,1)..(1+0.6,1);
\draw[thick] (0.4,0.1)..controls (0.45,0.10) and (0.75,0.435)..(0.8,0.5);
\draw[thick] (0.8,0.5)..controls (0.85,0.565) and (0.75,0.47)..(1.2,0.9);
\draw[thick] [->] (-1.5,0.0) -- (4,0.0);
\draw [thick][->] (0,-1.5) -- (0,2.5);
\draw[thick] (0+1.6,0+1)..controls (0.25+1.6,0+1) and (0.3+1.6,0.07+1)..(0.4+1.6,0.1+1);
\draw[thick](0.6+0.6+1.6,0.9+1)..controls (0.7+0.6+1.6,0.93+1) and (0.75+0.6+1.6,1+1)..(1+0.6+1.6,1+1);
\draw[thick] (0.4+1.6,0.1+1)..controls (0.45+1.6,0.10+1) and (0.75+1.6,0.435+1)..(0.8+1.6,0.5+1);
\draw[thick] (0.8+1.6,0.5+1)..controls (0.85+1.6,0.565+1) and (0.75+1.6,0.47+1)..(1.2+1.6,0.9+1);
\draw[thick] [->] (-1.5,0.0) -- (4,0.0);
\draw [thick][->] (0,-1.5) -- (0,2.5);
\draw [dotted, thick] (0.,1) -- (3.5,1);
\node [below] at (-0.31,1.1) {$\pi$};
\draw [dotted, thick] (0.,2) -- (3.5,2);
\node [below] at (-0.31,2.2) {$2\pi$};
\node [below] at (-0.2,0.1) {$0$};
\node [below] at (1.6,0.1) {$\frac12$};
\node [below] at (3.2,0.1) {$1$};
\node [below] at (4.2,0.2) {$x$};
\node [below] at (0,3.1) {$\xi(x)$};
\draw [dotted, thick] (1.6,0.) -- (1.6,1);
\draw [dotted, thick] (3.2,0.) -- (3.2,2);
\node[thick] [below] at (0.5,-1.5) {Figure 2: nonhomotopic to identity};
\end{tikzpicture}\ \ \
\end{center}

In the following, we will use $c,\  C,\  C(l)$, etc, to denote  universal positive constants independent of iterative steps.
For any cocycle $A(x)$,  $n\in Z^+$  and $x\in I$,
  we decompose $A^n(x)$ as $R_{\psi_{A,n}(x)}\cdot \Lambda_{A,n}(x)\cdot R_{\phi_{A,n}(x)}$  when $A^{n}(x)$ is hyperbolic in $I$ and decompose $A^n(T^{-n}x)$ as $R_{\psi_{A,-n}(x)}\cdot \Lambda_{A,-n}(x)\cdot R_{\phi_{A,-n}(x)}$ when $A^{n}(T^{-n}x)$ is hyperbolic in $I$.

Let $\xi_N(x)=\xi(x)$ defined above. Define $B_N(x)=\Lambda R_{\frac{\pi}{2}-\xi_N(x)}$.

  \begin{Proposition}\label{WY2}  There are ${\cal C}^{l}$ functions $\xi_k(x)$
($k=N+1, N+2, \cdots$) constructed inductively such that

 \begin{equation}\label{phin}\hskip -6cm {\bf{\it 1.}}\ |\xi_k(x)-\xi_{k-1}(x)|_{C^{l}}\le C(l)\cdot \lambda_k^{-2r_k}
\cdot |I_k|^{-l^2}.\end{equation}

2.\ Let $B_k(x)=\Lambda R_{\frac{\pi}{2}-\xi_k(x)}$.
For each $x\in I_k$, we have \beq\label{difference-potential}\|B_k^{ r^{\pm}_k(x)}(x)\|\ge \lambda_k^{r^{\pm}_k(x)}.\eeq
\ \ \ \ \ 3.\
For $x\in I_k$, we have
$$
\begin{array}{ll}
&{\rm (1)_k}\quad \psi_{B_k,-r^-_{k}}(x)+\phi_{B_k,r^+_{k}}(x)-\frac\pi2=\xi_0(x)\quad {\rm\ on}\ \frac{I_k}{10};\\
 \\
  &{\rm (2)_k}\quad |\psi_{B_k,-r^-_{k}}(x)+\phi_{B_k,r^+_{k}}(x)-\frac\pi2|\ge \frac{1}{(20q_{k}^2)^{l+1}},\quad
x\in I_k\backslash \frac{I_k}{10},
 \end{array}
$$
where $\xi_0(x)$ is defined in (\ref{phi0}) and (\ref{phi01}).
\end{Proposition}

\begin{Remark}\label{Remark2.1} It is easy to see from (\ref{phin}) that  $B_k$ converges to a limit $D_l$ in $\mathcal{C}^l$  -topology.
Moreover, from (\ref{phin}) and (\ref{difference-potential}) as well as Theorem 3 in \cite{WY}, we have $L(D_l)\ge (1-\epsilon)\ln \lambda$.\end{Remark}

\vskip 0.2cm

To prove Proposition \ref{WY2}, we first give the following Lemma \ref{iterative-estimate-derivative}.


\begin{Lemma}\label{iterative-estimate-derivative}
 For any function $\sigma(x)$ defined on $S^1$, let $d_k(\sigma)=\min_{x\not\in I_k}\{|\sigma(x)|\}.$ Assume that for any $x\in I_k$,
\beq\label{norm-n}\log \|A^{r_k}(x)\|\gg -\log d_{k+1},\eeq
where $d_{k+1}=d_{k+1}(\phi_{ A,r^+_{k}}(x)+\psi_{A,-r^-_{k}}(x)-\frac\pi2)$. Furthermore assume that, for $i\le l$ and $m^{\pm}=r^{\pm}_k(x)$,
$$\left\{\begin{array}{ll}
&\left|\frac{d^i}{dx^i}\phi_{A,m^+}(x)\right|, \quad\left|\frac{d^i}{dx^i}\psi_{A, -m^-}(x) \right|  \le C(i)\cdot d_{k+1}^{-i}     \qquad {\rm (1)}_k\\
 \\
&\left|\frac{d^i \|A^{\pm m}(x)\|}{dx^i}\right|\cdot\|A^{\pm m}(x)\|^{-1}\le C(i)\cdot d_{k+1}^{-i}.\qquad \quad{\rm (2)}_k
\end{array}\right.
$$
Then for $i\le l$, $x\in I_{k+1}$ and $\hat{m}^{\pm}=r^{\pm}_{k+1}(x)$
it holds that
$$\left\{\begin{array}{ll}
 &\left|\frac{d^i}{dx^i}\phi_{A,\hat{m}^+}(x)\right|, \left|\frac{d^i}{dx^i}\psi_{A, -\hat{m}^-}(x) \right|  \le C(i)\cdot d_{k+1}^{-i}, \qquad {\rm (1)}_{k+1}\\
 \\
&\left|\frac{d^i \|A^{\pm \hat m}(x)\|}{dx^i}\right|\cdot\|A^{\pm \hat m}(x)\|^{-1}\le C(i)\cdot d_{k+1}^{-i}.\qquad \quad{\rm (2)}_{k+1}
\end{array}\right.
$$

Moreover, for any $i\ge 0, x\in  I_{k+1}$, it holds that
\beq\label{difference-derivative}
\begin{array}{ll}
\left|\frac{d^i}{dx^i}(\phi_{A,r^+_{k+1}}(x)-\phi_{A,r^+_k}(x))\right|\le C(i)\cdot \|A^{r^+_k}\|^{-2}\cdot d_k^{-i},\\
\left|\frac{d^i}{dx^i}(\psi_{A,-r^-_{k+1}}(x)-\psi_{A,-r^-_{k}}(x))\right|\le C(i)\cdot \|A^{r^-_k}\|^{-2}\cdot d_k^{-i}.
\end{array}\eeq
\end{Lemma}
\vskip 0.3cm
The proof of Lemma \ref{iterative-estimate-derivative} will be given in the Appendix.
\\

\noindent
{\it Proof of Proposition \ref{WY2}.}
\ For each $k\ge N$ and $x\in I_k$, since
\[{\hat f}_k(x):=(\psi_{B_{k-1},-r^-_{k}}(x)+\phi_{B_{k-1},r^+_{k}}(x))-(\psi_{B_{k-1},-r^-_{k-1}}(x)+\phi_{B_{k-1},r^+_{k-1}}(x))\] usually does not vanish on $I_{k-1}$ and thus $\psi_{B_{k},-r^-_{k}}(x)+\phi_{B_{k},r^+_{k}}(x)-\frac\pi2\not= \xi_0(x)$ on $I_{k}$.  To guarantee ${\rm (1)_k}$ in Proposition \ref{WY2}, we modify $\xi_{k-1}(x)$ on $I_k$   as ${\xi}_{k}(x)=\xi_{k-1}(x)+f_k(x)$, where $C^l$ periodic function $f_k(x)$ is defined as follows
$$
f_{k}(x)=\left\{\begin{array}{ll}
{\hat f}_k(x) &\quad\  x\in \frac{I_{k}}{10}\\
\\
h^{\pm}_{k}(x), &\quad\  x\in I_{k}\backslash \frac{I_{k}}{10}\\
\\
0,& x\in \mathbb{S}^1\backslash I_{k}
\end{array}\right.
$$
where $h^{\pm}_{k}(x)$ is a polynomial of degree $2l+1$ restricted in each interval of
$I_{k}\backslash \frac{I_{k}}{10}$ satisfying
$$\begin{array}{ll}
&\frac{d^jh^{\pm}_{k}}{dx^j}(\pm\frac{1}{10q_{k}^2})=\frac{d^j{\hat f}_k}{dx^j}(\pm\frac{1}{10q_{k}^2})\\
\\
&\frac{d^jh^{\pm}_{k}}{dx^j}(\pm\frac{1}{q_{k}^2})=0,\quad i=1,
2, \quad 0\le j\le l.
\end{array}
$$
From (\ref{difference-derivative}) in Lemma \ref{iterative-estimate-derivative}, we have
\beq\label{hatf_k}
|(\psi_{B_{k-1},-r^-_{k}}(x)+\phi_{B_{k-1},r^+_{k}}(x))-(\psi_{B_{k-1},-r^-_{k-1}}(x)+\phi_{B_{k-1},r^+_{k-1}}(x))|_{C^l}\le C(l)\cdot \lambda_k^{-2r_k}
\cdot |I_k|^{-l^2},
\eeq
where (\ref{norm-n}) is fulfilled by conclusion 2 and 3 of the induction assumption for the case $k-1$.

In view of  the definition of $f_k(x)$ we obtain
\beq\label{f_k}
|f_k|_{C^l}\le C(l)\cdot \lambda_k^{-2r_k}
\cdot |I_k|^{-l^2}.
\eeq

Let $B_{k}(x)=\Lambda \cdot R_{\frac{\pi}{2}-\xi_{k}(x)}$, then we have

\begin{Lemma}\label{rotationn}
For $x\in I_k$, it holds that
$$
B_k^{ {r_k^+(x)}}(x)=B_{k-1}^{ {r_k^+(x)}}(x)\cdot R_{-f_k(x)}
 $$
and
$$
B_k^{r_k^-(x)}(T^{-r_k^-(x)}x)= B_{k-1}^{r_k^-(x)}(T^{-r_k^-(x)}x).
 $$
\end{Lemma}
\Proof Obviously $T^ix\in \mathbb{S}^1\backslash I_k$ for $x\in I_k$ and
$1\le i\le  {r_k^+}(x)-1$. Since $B_k(x)=B_{k-1}(x)$ for $x\in
\mathbb{S}^1\backslash I_k$, we have that
$$
B_k^{ {r_k^+(x)}}(x)=B_{k-1}^{ {r_k^+(x)}}(x)\cdot (B_{k-1}^{-1}(x)B_k(x)),\quad
x\in I_k.
$$
From the definition, we have $B_k(x)=B_{k-1}(x)\cdot
R_{\xi_{k-1}(x)-\xi_k(x)}$, which implies
$B_{k-1}^{-1}(x)B_k(x)=R_{\xi_{k-1}(x)-\xi_k(x)}$. Thus we obtain the first
equation in Lemma \ref{rotationn}. Similarly we can prove the second one.\hfill{}
\qedbox
\vskip 0.3cm
\begin{Lemma}\label{tophi}
It holds that
$$f_k(x)=(\psi_{B_{k-1},-r^-_{k}}(x)+\phi_{B_{k-1},r^+_{k}}(x))-(\psi_{B_{k},-r^-_{k}}(x)+\phi_{B_{k},r^+_{k}}(x)),\quad x\in I_k.$$
\end{Lemma}
\Proof Since a rotation does not change the norm of a vector, for a
hyperbolic matrix $A$ and a rotation matrix $R_{\t}$, we have
\beq\label{factonrotation} \phi_{A\cdot R_{\t}}=\phi_{A}+\t. \eeq From Lemma \ref{rotationn}, we have
$$
\phi_{B_{k},r^+_{k}}(x)=\phi_{B_{k-1},r^+_{k}}(x)-{f_k(x)}, \quad \psi_{B_{k},-r^-_{k}}(x)=\psi_{B_{k-1},-r^-_{k}}(x).
$$
Thus
$$f_k(x)=(\psi_{B_{k-1},-r^-_{k}}(x)+\phi_{B_{k-1},r^+_{k}}(x))-(\psi_{B_{k},-r^-_{k}}(x)+\phi_{B_{k},r^+_{k}}(x)),\quad
x\in I_k,$$ which concludes the proof.\hfill{}  \qedbox \vskip 0.6cm

\noindent
{\it  Proof of $(1)_k$ and $(2)_k$ }\quad  From the definition of $f_k(x)$, we have
$f_k(x)=(\psi_{B_{k-1},-r^-_{k}}(x)+\phi_{B_{k-1},r^+_{k}}(x))-(\psi_{B_{k-1},-r^-_{k-1}}(x)+\phi_{B_{k-1},r^+_{k-1}}(x))$ on $\frac{I_k}{10}$, which together with
Lemma \ref{tophi} implies that for each $x\in \frac{I_k}{10}$,
$$\psi_{B_{k},-r^-_{k}}(x)+\phi_{B_{k},r^+_{k}}(x)=(\psi_{B_{k-1},-r^-_{k}}(x)+\phi_{B_{k-1},r^+_{k}}(x))-f_k(x)=
\psi_{B_{k-1},-r^-_{k-1}}(x)+\phi_{B_{k-1},r^+_{k-1}}(x).$$ Since
$\psi_{B_{k-1},-r^-_{k-1}}(x)+\phi_{B_{k-1},r^+_{k-1}}(x)=\xi_0(x)$ on $\frac{I_{k-1}}{10}$ by induction assumption $(1)_{k-1}$, we obtain $(1)_k$ in proposition \ref{WY2}.

 Obviously $\lambda_k^{q_{k-1}}\gg q_k^{2l}$.
Hence  ${\rm (2)_k}$ in Proposition \ref{WY2} can be obtained from the
induction assumption $(2)_{k-1}$ and (\ref{f_k}).
\vskip 0.3cm

\noindent
{\it Proof of conclusion 1 of Proposition \ref{WY2}.}\quad
Conclusion 1 can be obtained from  (\ref{hatf_k}).
 \vskip 0.3cm
 \noindent
 {\it Proof of conclusion 2 of Proposition \ref{WY2}.}\quad   For $x\in I_k$, let $i_1(x)<i_2(x)<\cdots <i_{j(x)}(x)\le r_k$ be the returning times of $I_{k-1}$ less than $r_k$.
    Since $|I_{k}|\le \frac14|I_{k-1}|$ (we can make a slight modification of the definition of $I_{k}$ if necessary), from the symmetry between $I_{k,1}$ and $I_{k,2}$,  we have that for any $x\in I_{k}$, we have $T^{r_{k}}x\in I_{k-1}$. Then  we have that $i_{j(x)}(x)=r_k$. Since $T^{i_s(x)}x\not\in I_k$ for $s< {j(x)}$,  $|\theta_s-\frac\pi2|\ge \frac{1}{q_k^{2l}}$, where $\theta_s=\phi_{B_k, i_{s+1}(x)-i_{s}(x)}(T^{i_s(x)}x)+\psi_{B_k,i_s(x)-i_{s-1}(x)}(T^{i_{s-1}(x)}x)$. Together with the conclusion 3 of the induction assumption for $(k-1)$-th step we have that
    $|\tilde{\theta}_s-\frac\pi2|\ge \frac{1}{2q_k^{2l}}$, where $\tilde{\theta}_s=\phi_{B_k, i_{s+1}(x)-i_{s}(x)}(T^{i_s(x)}x)+\psi_{B_k,i_s(x)}(x).$
    Thus from the definition of $\l_k$,  we
 obtain the conclusion 2 for $k$-th step by  repeated applications of Lemma \ref{L2.1}.

\begin{Remark}
In spirit, the proof of conclusion 2 of Proposition \ref{WY2} coincides with the one of LDT.
\end{Remark}

\section{The construction of $C_k(x)$}\label{s4}

 Now we start to construct a $ C_k$ in any small $\mathcal{C}^l$-neighborhood of $B_k$ such  that $L(C_k)=0$. It is obvious that $C_k\rightarrow D_l$ in $\mathcal{C}^l$ topology. $C_k$ will be constructed as  limit of a sequence of converging cocycles, say $A_{k,i}$, in any small neighborhood of $B_k$ such that $L(A_{k,i})\to 0$ as $i\to\infty$. By the construction, we can show that $L(C_k)= \lim_{i\rightarrow \infty} L( A_{k,i})=0$, see Corollary \ref{Cor4.1}. In the following, we shall simply denote $A_{k, i}$ by $A_i$.

  \vskip 0.2cm
 The following lemma is of key importance for the construction:
\vskip 0.2cm
\noindent {\bf Iterative Lemma:}\ {\it Let $A_0(x)=\Lambda \cdot R_{\frac\pi2-\t_0(x)}$ satisfy that $\|A_0^{r_{n_0}(x)}(x)\|\ge \mu^{r_{n_0}(x)}$ with $\lambda\ge\mu\gg1$, $n_0\ge N$ and $\psi_{A_0, -r_{n_0}}(x)+\phi_{A_0, r_{n_0}}(x)-\frac\pi2=\xi_0(x)$, $x\in I_{n_0}$. Then we can find two small positive numbers ${\delta}_1>\delta_2$ such that for any $i\ge 0$, there exist $A_i(x)=\Lambda\cdot R_{\frac\pi2-\t_i(x)}$ and ${n_i}$, such that the following hold
$$
{\rm(P_i)}:\ \left\{\begin{array}{ll} {\rm (1).} & \|A_i^{r_{n_i}(x)}(x)\|\ \ge \mu^{(1-\delta_1)^i\cdot r_{n_i}(x) } {\rm \ on}\ I_{n_i}\ {\rm and}\ \mu^{(1-\delta_1)^i\cdot q_{n_i}}\gg \frac{1}{|I_{n_i}|};\\
{\rm (2).} &  \|A_i^{r_{n_j}(x)}(x)\|\le \lambda^{(1-\delta_2)^j\cdot{{r_{n_j}(x)}}} {\ \rm for\ }  x\in I_{n_j}\ {\rm and\  } j\le i;\\
{\rm (3).} & \bar{\mu}_{n_i}\le \underline{\mu}_{n_i}^{2};\\
{\rm (4).} & \psi_{A_i, -r_{n_i}}(x)+\phi_{A_i, r_{n_i}}(x)-\frac\pi2=\xi_0(x)\ {\rm on}\ I_{n_i};\\
{\rm (5).} & |\t_{i+1}-\t_i|_{\mathcal{C}^l}\le q_{n_i}^{4Ml^2}\cdot \mu^{-\frac12(1-\delta_1)^i\cdot q_{n_i}}+q_{n_i}^{-2}.
\end{array}
\right.
$$
In the above, $\bar{\mu}_{n_i}=\max_{x\in I_{{n_i}}} \|A_i^{r_{n_i}(x)}(x)\|^{\frac{1}{r_{n_i}(x)}}$ and $\underline{\mu}_{n_i}=\min_{x\in I_{{n_i}}} \|A_i^{r_{n_i}(x)}(x)\|^{\frac{1}{r_{n_i}(x)}}.$
Therefore, $\underline{\mu}_{n_i}\ge \mu^{(1-\delta_1)^i}$ and $\bar{\mu}_{n_i}\le \lambda^{(1-\delta_2)^i}$.
\vskip 0.5cm

\noindent  {\rm The main result Theorem \ref{maintheorem} is an easy consequence of the following corollary.}
\noindent \begin{Corollary}\label{Cor4.1} There exists a $SL(2, \mathbb{R})$-sequence $\{C_k\}_{k=N}^{\infty}$ with $L(C_k)=0$ such that $C_k$ tend to $D_l$ in the $\mathcal{C}^l$ topology.
\end{Corollary}
\noindent \begin{Proof}\ For any $k\in \mathbb{N}$, we  apply Iterative Lemma by setting $A_0=B_k$, $n_0=q_k$ and $\mu=\lambda^{1-\epsilon}$ where $B_k$ is defined in Proposition \ref{WY2}. Hence for each $i$ we obtain $A_i$ such that $\rm {(P_i)}$\ holds true. By (5) and the inequality $\mu^{(1-\delta_1)^i\cdot q_{n_i}}\gg \frac{1}{|I_{n_i}|}$ in (1) of $\rm {(P_i)}$, $A_i$ has a limit, say $C_k$, in the ${\cal C}^l$ topology. From (2) of $\rm {(P_i)}$, we obtain $\|C_k^{r_{n_j}(x)}(x)\|^{\frac{1}{r_{n_j}(x)}}\le \l^{(1-\delta_2)^j}$ for any $j\le i$ and  $x\in I_{n_j}$, which by the subadditivity of Lyapunov exponent, implies $L(C_k)\le (1-\delta_2)^j\log\l$ for any $j$. Let $j\rightarrow \infty$, we obtain
   $L(C_k)=0.$  Moreover, from (5) of  $\rm {(P_i)}$ it holds that
$$
\begin{array}{ll}
\|C_k-D_l\|_{\mathcal{C}^l}&\le \|C_k-B_k\|_{\mathcal{C}^l}+\|B_k-D_l\|_{\mathcal{C}^l}\\
\\
&=\|C_k-A_{0}\|_{\mathcal{C}^l}+\|B_k-D_l\|_{\mathcal{C}^l}\\
\\
&\le 2(q_k^{4M^2l}\cdot\l^{-(1-\epsilon)(1-\delta_1)\cdot q_k}+q_k^{-2})+\|B_k-D_l\|_{\mathcal{C}^l}.
\end{array}
$$
In the last inequality, we use (5) and the inequality $\mu^{(1-\delta_1)^i\cdot q_{n_i}}\gg \frac{1}{|I_{n_i}|}$ in (1) of $\rm {(P_i)}$. It implies $C_k\rightarrow D_l$ in the $\mathcal{C}^l$ topology as $k\rightarrow \infty$.
Moreover,   $
L(D_l)\ge (1-\epsilon)\ln \lambda$ by Remark \ref{Remark2.1}.\qquad \hfill{}\qedbox
\end{Proof}\\

\noindent
{\bf Proof of Iterative Lemma.}
{\rm Let $0<\delta_0\ll \max\{\frac{1}{l^2}, M^{-k_1}\}$ be a fixed number and $\delta_1=8 \delta_0l$, $\delta_2=M^{-k_1}\cdot \delta_0l$, where $k_1$ is defined in Proposition \ref{ratio}.

When $i=1$, }${\rm (P_i)}$ {\rm obviously holds true for $A_1$ with $\lambda\gg 1$.
 Assuming that $A_1, \cdots, A_{i-1}$ have been constructed with $\rm (P_1),\cdots, (P_{i-1})$, we will construct $A_i$ such that $\rm (P_i)$ holds.  From (3) of ($P_{i-1}$), we have $ \|A^{r_{n_{i-1}}}_{i-1}(x)\|\le \|A^{r_{n_{i-1}}}_{i-1}(y)\|^{{2}}$ for $x,\ y\in I_{n_{i-1}}$.

\noindent {\it Step 1. Definition of $n_{i}$ and $I_{n_i}$}.\ \ \ \
{\rm
 Choose  $n_i\gg n_{i-1}$  such that  $\mu^{(1-\delta_1)^{i}\cdot q_{n_i}}\gg q_{n_i}^{2}\gg \lambda^{2\delta_0\cdot {r_{n_{i-1}}}}$ and $I_{n_i}$ is defined as before. The Diophantine condition implies that $r_{n_i}\ge q_{n_i}$.


\vskip 0.3cm
\noindent{\it Step 2. Modification of $A_{i-1}$.}\ \ \
For our purpose, we first make a local modification for $A_{i-1}$ on $I_{n_{i-1}}$  such that there is a low platform in the image of $\phi_{ \tilde{A}_{i-1},r_{n_{i-1}}}(x)+\psi_{\tilde{A}_{i-1},-r_{n_{i-1}}}(x)-\frac\pi2$ for the new cocycle $\tilde{A}_{i-1}$, see Figure 3.


Consider the sub-interval $[0, \frac{1}{{q_{n_{i-1}}^2}}]$ of $I_{n_{i-1}}$. Define $0<c
<\tilde{c}< d<\frac{1}{{q_{n_{i-1}}^2}}$ such that
$|c|={\underline\mu}_{n_{i-1}}^{-2\delta_0\cdot {r_{n_{i-1}}}}, |\tilde{c}|=M^2\cdot |c|, {d}=\frac{1}{2{q_{n_{i-1}}^2}}$, see Figure 3. From the definition of $n_i$, we have $|I_{n_i}|< |c|.$

 Define
$$
e_i^0(x)=\left\{\begin{array}{ll}
2|c|^{l+1}-(\phi_{ {A}_{i-1},r_{n_{i-1}}}(x)+\psi_{{A}_{i-1},-r_{n_{i-1}}}(x)-\frac\pi2),& x\in [\tilde{c}, {d}];\\
0, & x\notin [c, \frac{1}{{q_{n_{i-1}}^2}}]\\
\tilde{h}_i(x) & x\in [c, \tilde{c}]\cup\ [d, \frac{1}{{q_{n_{i-1}}^2}}],
\end{array}
\right.
$$
where $\tilde{h}_{i}(x)$ are  polynomials of degree $2l+1$
restricted on each interval
and for $0\le j\le l$ satisfies
$$\begin{array}{ll}
\frac{d^j\tilde{h}_{i}}{dx^j}(
\frac{1}{{q_{n_{i-1}}^2}})=0,&\frac{d^j\tilde{h}_{i}}{dx^j}(d)=\frac{d^j(2|c|^{l+1}-(\phi_{ {A}_{i-1},r_{n_{i-1}}}+\psi_{{A}_{i-1},-r_{n_{i-1}}}-\frac\pi2))}{dx^j}(d),\\
\\ \frac{d^j\tilde{h}_{i}}{dx^j}(
c)=0,&\frac{d^j\tilde{h}_{i}}{dx^j}(\tilde{c})=\frac{d^j(2|c|^{l+1}-(\phi_{ {A}_{i-1},r_{n_{i-1}}}+\psi_{{A}_{i-1},-r_{n_{i-1}}}-\frac\pi2))}{dx^j}(\tilde{c}).\end{array}
$$
$e_i^1(x)$ on the subinterval $[\frac 12, \frac 12+\frac{1}{{q_{n_{i-1}}^2}}]$ of $I_{n_{i-1}}$ is defined  similarly.
Let
$e_i(x)=e_i^0(x)+e_i^1(x)$.
We have the following estimates for $e_i(x)$.
\begin{Lemma}\label{L4.1}
It holds that $|{e}_{i}(x)|_{{\cal
C}^l}\le C\cdot{q_{n_{i-1}}^{-2}}$.
\end{Lemma}
\begin{Proof}
From (1) in $\rm (P_{i-1})$ we have $\mu_{n_{i-1}}^{r_{n_{i-1}}}\ge \mu^{(1-\delta_1)^i\cdot q_{n_{i-1}}}\gg \frac{1}{|I_{n_{i-1}}|}$. Then From (4) in $\rm (P_{i-1})$ and the definition of $c$, it holds for $0\le j\le l$ that
$$|(2|c|^{l+1}-(\phi_{ {A}_{i-1},r_{n_{i-1}}}+\psi_{{A}_{i-1},-r_{n_{i-1}}}-\frac\pi2))(x)|_{{\cal C}^j}\le C\cdot{q_{n_{i-1}}^{-2(l+1-j)}}.$$ Hence from
Cramer's rule we have that $|\tilde{h}_{i}(x)|_{{\cal
C}^l}\le C\cdot{q_{n_{i-1}}^{-2}}$. Consequently, $|{e}_{i}(x)|_{{\cal
C}^l}\le C\cdot {q_{n_{i-1}}^{-2}}$.
\end{Proof}\hfill{}
\qedbox

Let $\tilde{\t}_i=\t_{i-1}+e_i(x)$ and  $\tilde{A}_{i-1}=\Lambda\cdot
R_{\frac\pi2-\tilde{\t}_i}$, we have $\psi_{\tilde{A}_{i-1},-r_{n_{i-1}}}(x)+\phi_{\tilde{A}_{i-1},r_{n_{i-1}}}(x)-\frac\pi2$ on the part $[0, \frac{1}{q^2_{n_{i-1}}}]$ of $I_{n_{i-1}}$ is of the shape in Figure 3.

\begin{figure}
 \centering
 \resizebox{6.7in}{2.9in}
{ \includegraphics{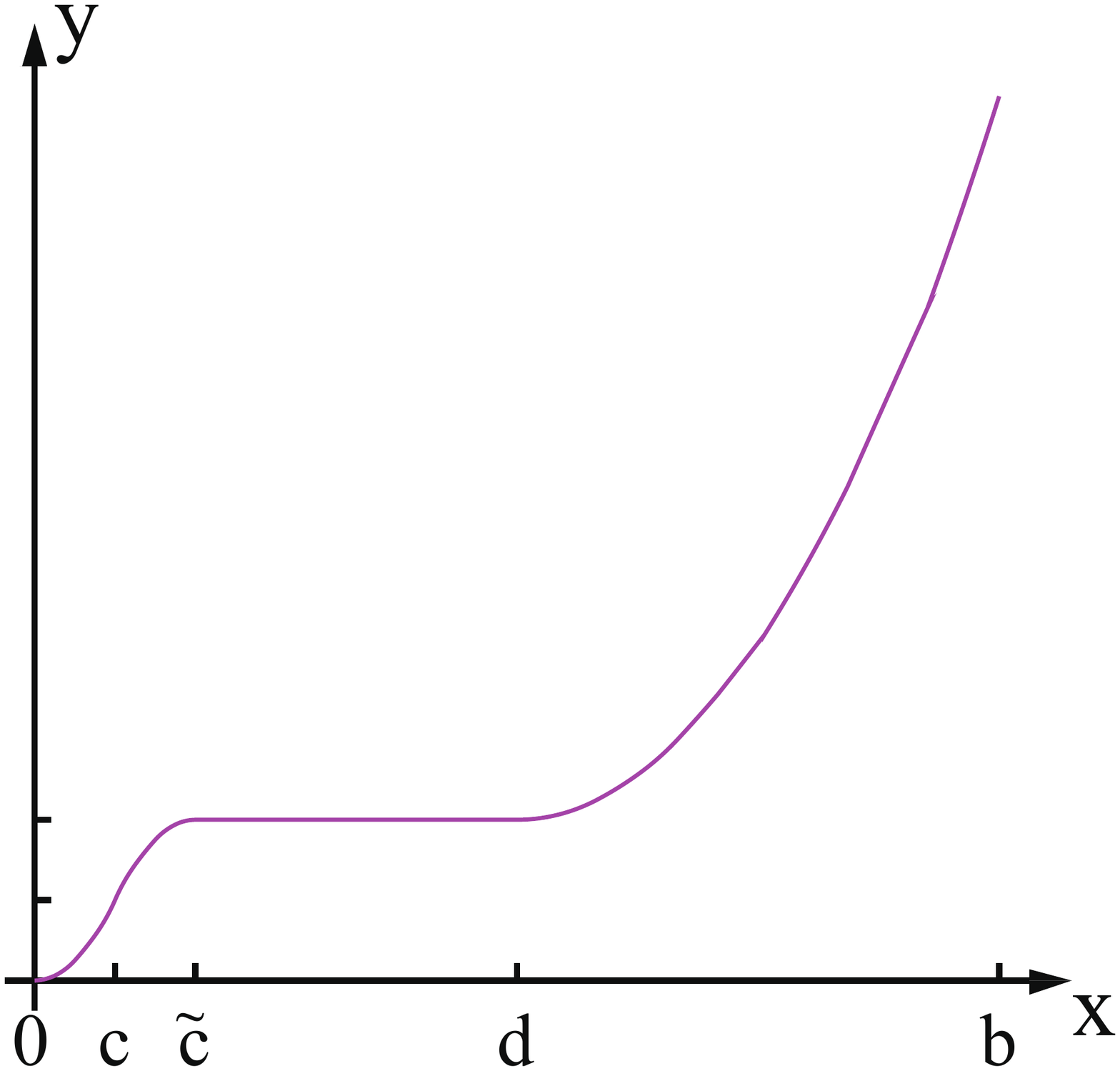}}\begin{center}
Figure 3, $b=\frac{1}{{q_{n_{i-1}}^2}}$\end{center}
\end{figure}


\vskip 0.5cm
\noindent{\it Step 3. The estimate on the lower bound.}\ \ \
\begin{Lemma}\label{sub-iterative} Let $A_{i,0}=\Lambda \cdot R_{\tilde{\t}_i}:=\Lambda \cdot R_{\t_{i,0}}$ satisfy
$\| A_{i,0}^{r_{n_{i-1}}(x)}(x)\|\ge\nu_{0}^{r_{n_{i-1}}(x)}$  for $x\in I_{n_{i}}$ with $\nu_0=\mu^{(1-\delta_1)^{i-1}}$. Then for any $n_i-n_{i-1}\ge j\ge 1$,
there exist $\t_{i,j}$ and $A_{i,j}=\Lambda\cdot R_{\frac\pi2-\t_{i,j}}$ such that the following properties hold true :

$$
{\rm(\tilde{P}_{i,j})}:\ \left\{\begin{array}{ll} {\rm \tilde{(1)}.} &
\| A_{i,j}^{r_{n_{i-1}+j}(x)}(x)\|\ge \nu_{j}^{r_{n_{i-1}+j}(x)} {\rm on}\ I_{n_i};
\\
 {\rm \tilde{(2)}.} & \phi_{A_{i,j},r_{n_{i-1}+j}}(x)+\psi_{A_{i,j}, -r_{n_{i-1}+j}}(x)-\frac\pi2=\t_{i,0}(x)\ {\rm on}\ I_{n_{i}};\\
{\rm \tilde{(3)}.} & |\t_{i,j}-\t_{i, j-1}|_{\mathcal{C}^l}\le r_{i,j}\cdot \nu_{j-1}^{-q_{n_{i-1}+j-1}},\ r_{i,j}\approx \max\{\nu_{j-1}^{2l^2\delta_0q_{n_{i-1}}}, q^{2l^2}_{n_{i-1}+j}\},
\end{array}
\right.
$$
where  $\nu_j$ are iteratively defined by
$$\nu_{j}=\nu_0\cdot\nu_0^{-\delta_0
(\sqrt{2}^{-(j-1)}+\cdots +\sqrt{2}^{-1}+1)\cdot 2(l+1)}\ge \nu_0^{(1-8\delta_0\cdot (l+1))}= \mu^{(1-\delta_1)^i}.$$
Let $\underline{\mu}_{i,j}=\min_{x\in I_{{n_i}}} \|({A}^{r_{n_{i-1}+j}}_{i,j}(x))\|^{\frac{1}{{r_{n_{i-1}+j}}}}$ for any $j\le n_i-n_{i-1}$. We have $\underline{\mu}_{i,j}\ge \nu_{j}$ and $\underline{\mu}_{n_i}=\underline{\mu}_{i,n_i-n_{i-1}}\ge  \nu_{n_i-n_{i-1}}$.


\end{Lemma}

\noindent \begin{Proof}\quad
For $j=1$, from (1) of ($P_{i-1}$) and the definitions of $\tilde{\t}_i$ and $\mu_0$, we have
$$\begin{array}{ll}
&\frac{1}{r_{n_{i-1}+1}(x)}\log\|A_{i,1}^{r_{n_{i-1}+1}(x)}(x)\|\ge \log{\underline\mu}_{i,n_{i-1}}+ \frac{1}{q_{n_{i-1}+1}} \log{\underline\mu}_{i,n_{i-1}}^{-2(l+1)\delta_0 q_{n_{i-1}}\cdot \frac{q_{n_{i-1}+1}}{q_{n_{i-1}}}}\\
\\
&=(1-2(l+1)\delta_0)\log{\underline\mu}_{i,n_{i-1}} \ge \nu_1.\end{array}$$
Thus we obtain $\tilde{(1)}$. Moreover $\tilde{(2)}$ and $\tilde{(3)}$ can be proved by Proposition \ref{WY2} and Lemma \ref{iterative-estimate-derivative} with $d_{n_{i-1}+1}\ge \frac{1}{q^{2(l+1)}_{n_{i-1}+1}}$.

Assume $(\tilde{P}_{i,j})$ hold true. We will prove $(\tilde{P}_{i,j+1})$. Define $\t_{i,j+1}(x)$ by modifying  $\t_{i,j}(x)$ in the same way as we define $\xi_{k+1}$ by modifying $\xi_k$ in Proposition \ref{WY2}. Applying Lemma \ref{iterative-estimate-derivative} with $d_{n_{i-1}+j}\ge \min\{\nu_{j-1}^{-2(l+1)\delta_0q_{n_{i-1}}}, q^{-2(l+1)}_{n_{i-1}+j}\}$, we get $\tilde{(2)}$ and $\tilde{(3)}$.

Now we prove  $\tilde{(1)}$. In case that $|\tilde{c}|<|I_{n_{i-1}+j}|$, we have $q_{m+1}\ge \sqrt{2} \cdot q_m$ for each $m$, and thus
$$\begin{array}{ll}
&\frac{1}{r_{n_{i-1}+j}(x)}\log\|A_{i, j+1}^{r_{n_{i-1}+j}(x)}(x)\|\ge \log\nu_j+\frac{1}{q_{n_{i-1}+j+1}}\cdot \log{\underline\mu}_{i,n_{i-1}}^{-\delta_0\cdot q_{n_{i-1}}
\cdot \frac{q_{n_{i-1}+j+1}}{q_{n_{i-1}+j}}\cdot 2(l+1)}\\
\\
&\ge ({1-\delta_0
(\sqrt{2}^{-(j-1)}+\cdots +\sqrt{2}^{-1}+1)\cdot 2(l+1)})\log\nu_0{-\delta_0 \cdot \sqrt{2}^{-j}\cdot 2(l+1)}\log\nu_0\\
\\
&=\log\nu_{j+1}.
\end{array}
$$
Now we consider the case $|\tilde{c}|\ge |I_{n_{i-1}+j}|$.
Let $j^*$ be the smallest integer such that $|I_{q_{n_{i-1}+j^*}}|\le |\tilde{c}|$ (Obviously, $j^*$ depends on $n
_{i-1}$ and we can choose $n_i$ large enough such that $j^*\ll n_i$). Since for any $s$, it holds that
$M^2\cdot |I_{s+1}|\ge |I_s|$. Thus from the definition of $|c|$ and $|\tilde{c}|$, we have $|I_{q_{n_{i-1}+j^*}}|\ge |c|$ since $|I_{q_{n_{i-1}+j^*-1}}|\ge |\tilde{c}|$.
 Notice that $\nu_{j^*}^{q_{n_{i-1}+j^*}}\gg q_{n_{i-1}+j^*+1}^{2l}$. We construct
$\psi_{i, j^*+m}$ and $A_{i, j^*+m}=\Lambda \cdot R_{\psi_{i, j^*+m}}$  as in Proposition \ref{WY2},such that
$$
\underline{\mu}_{i,n_{i-1}+j^*+m}\ge \nu_{j^*+m}, \quad {\rm for}\ \  m\ge 1, $$
 which thus implies $\tilde{(1)}$.  \end{Proof}\qedbox
\vskip 0.2cm
\noindent Define
$\t_i(x)=\t_{i, j^*+m^*}(x)$ and $A_i(x)=\Lambda \cdot R_{\frac\pi2-\t_i(x)}$, where $m^*=n_i-n_{i-1}-j^*$. Then
(1) of ${\rm (P_i)}$\ can be proved by  $\tilde{(1)}$ in
$(\tilde{P}_{i,j})$. From the inequality $0<\delta_0\ll \frac{1}{l^2}$, (5) of ${\rm (P_i)}$\ can be proved by  $\tilde{(3)}$ in
$(\tilde{P}_{i,j})$ and Lemma \ref{L4.1}. (4) of ${\rm (P_i)}$\ is obvious from the construction of $A_{i,j}$.
\vskip 0.6cm

\noindent{\it Step 4.  The estimate on the upper bound.}\quad

\noindent Now we  prove (2) of ${\rm (P_i)}$, i.e., an upper bound estimate for the Lyapunov exponent.
We need the following proposition in \cite{WY}:
\begin{Proposition}\label{ratio}
Let $I_1$ be a small interval in $S^1$, $I_2=I_1+1/2$, $I=I_1\bigcup I_2$. Let $$\min r(x)=\min_{x\in I} \min \{i>0|  T^i x\ (mod\ 2\pi) \in
I\},$$  $$\max {r}(x)=\max_{x\in \frac 1{10} I_{1}} \min \{i>0|  T^i
x\ (mod\ 2\pi) \in \frac 1 {10}{I}_{1}\}.$$ Then there exists $k_1\in \mathbb{N}$ such that $M^{-k_1} \le
\frac{\min r(x)}{\max {r}(x)}\le 1$.
\end{Proposition}

Obviously,
we have  \beq\label{shunshi}|\phi_{A_i, r_{n_{i-1}}}(x)+\psi_{A_i, -r_{n_{i-1}}}(x)-\frac\pi2|\le 2|c|^{l}=2{\underline\mu}_{n_{i-1}}^{-2l\delta_0q_{n_{i-1}}}\eeq for $x\in [0,\frac{1}{2{q_{n_{i-1}}^2}}]$.
 Apply Proposition \ref{ratio} with $I_1=[0, \frac{1}{2{q_{n_{i-1}}^2}}]$. Then from (\ref{shunshi}) and Lemma \ref{L2.1} we have that
  $$\bar{\mu}_{n_{i}}\le \bar{\mu}_{n_{i-1}}\cdot\underline{\mu}_{n_{i-1}}^{-2lM^{-k_1}\cdot \delta_0}.$$
  Subsequently (3) of ($P_{i-1}$) and the definition of ${\delta}_2$ imply that
 $$\bar{\mu}_{n_{i}}\le \bar{\mu}_{n_{i-1}}^{1-\delta_2}\le \lambda^{(1-\delta_2)^i}.$$

\noindent{\it Step 5. The comparison between the lower and upper bounds.} \ \ For (3) of ($P_0$), the upper bound for  $w_0$ can be achieved by choosing $\lambda\gg N\gg 1$. For $i>0$, we have the following argument. For any $x, y$ with $|x-y|\le \frac{1}{q_{n_{i}}^2}$ we have that $\|A_i(x)-A_i(y)\|_{{\mathcal C}^l}\le \frac{C}{q_{n_{i}}^2}$. From (1) of ($P_{i-1}$) it holds that $\mu^{(1-\delta_1)^i\cdot q_{n_{i-1}}}\gg \frac{1}{|I_{n_{i-1}}|}$.
Once $n_{i-1}$ is determined, for any $\tilde{\epsilon}>0$, we can find $n_i\gg\tilde{n}_i\gg n_{i-1}$ such that
 $|I_{\tilde{n}_i}|\ll |c|$ and for any $x, y$ with $|x-y|\le \frac{1}{q_{{n}_{i}}^2}$, it holds that
$$\begin{array}{cc}&({1-\tilde{\epsilon}})\cdot L_{r_{\tilde{n}_i}}(A_i(y))\le L_{r_{\tilde{n}_i}}(A_i(x))\le ({1+\tilde{\epsilon}})\cdot L_{r_{\tilde{n}_i}}(A_i(y)),\\
\\
&\|\phi_{A_i, \pm r^{\pm}_{\tilde{n}_i}}(x)-\phi_{A_i,\pm r^{\pm}_{\tilde{n}_i}}(y)\|_{\mathcal{C}^l}+\|\psi_{A_i, \pm r^{\pm}_{\tilde{n}_i}}(x)-\psi_{A_i,\pm r^{\pm}_{\tilde{n}_i}}(y)\|_{\mathcal{C}^l}\le \tilde{\epsilon}.
\end{array}$$

Apply the inductive process for $A_i$ from step $\tilde{n}_{i}$ to $n_i$ for $x\in I_{n_i}$.  Thus similar to Remark \ref{Remark2.1}, if $\tilde{n}_i$ is large enough such that $\mu^{(1-\delta_1)^i\cdot q_{\tilde{n}_i}}\gg \frac{1}{|I_{\tilde{n}_i}|}$, we have that for any $x\in I_{n_i}$,
$$
\prod_{j} \|A^{{r_{\tilde{n}_{i,j}}}}_i(x_j)\|\ge \|A^{r_{n_i}}_i(x)\|\ge (\prod_{j} \|A^{{r_{\tilde{n}_{i,j}}}}_i(x_j)\|)^{1-\tilde{\epsilon}},
$$
where $x_j$'s are the points on orbits of $x$ when returning to $I_{\tilde{n}_i}$ and $\tilde{n}_{i,j}$'s are the corresponding returning times.
Hence for any $x, y\in I_{n_i}$ we have
$$
(1-2\tilde{\epsilon}) L_{r_{n_i}}(A_i(x))\le L_{r_{{n}_i}}(A_i(y))\le (1+2\tilde{\epsilon}) L_{r_{n_i}}(A_i(x)).
$$
 Then it follows that (3) of ($P_i$) holds true. This ends the proof of Iterative Lemma.}
\hfill{} \qedbox


\section{The proof for the ${\cal C}^{\infty}$ case}\label{s5}
{\rm
In this section, we will prove Theorem \ref{T2} and \ref{maintheorem} for the ${\cal C}^{\infty}$ case.
The basic idea is same as the one for the finitely smooth case. Essentially, we only need to modify cocycles in $C^\infty$ category.
 We will focus on the difference between the two cases.
 First we  follow the steps in section 3 to construct a sequence of $C^\infty$ cocycle  which is ${\cal C}^1$-convergent. Then we will prove that it actually converges  in ${\cal C}^{\infty}$ topology.

Assume $\l\gg e^{q_N^{a+1}}\gg 1$ with $0<a<\frac{1}{10}$. For
$n>N$, let $\l_{n+1}^{q_{n+1}}=\l_{n}^{q_{n+1}}\cdot
e^{-(10q_{n+1}^2)^{a}}$ with $\l_N=\l$.  From the definition of
$\l_n$, we have $\l_n^{q_n}\ge \l_{n-1}^{q_n}\cdot e^{-q_n^{2a}}\ge
\l_{n-2}^{q_n}\cdot
 e^{-q_n\cdot q_{n-1}^{2a-1}}\ge \cdots \ge \l^{q_n}\cdot\l_N^{-C\cdot q_n^{2a}}\ge \l^{(1-\e)q_n}$ for some small positive $\e$ if $\l\gg 1$ and $N\gg 1$. It implies that $\l_n$ decrease to $\l_{\infty}>\l^{1-\e}$.

 \noindent {\it Construction of $B_N(x)$}\qquad Let

\noindent(a)$$\xi_0(x)=\left\{\begin{array}{ll}
&\xi_{01}(x) \quad {\rm for}\  |x|\le \delta,\\
&\xi_{02}(x) ({\rm or\ } -\xi_{02}(x)) \quad {\rm for}\  |x-1/2|\le \delta,
\end{array}\right.
$$
where $\xi_{01}(x)={\rm sgn}(x)e^{-\frac{1}{|x|^a}}$ and $\xi_{02}(x)={\rm sgn}(x-1/2)e^{-\frac{1}{|x-1/2|^a}}$, $\delta>0$ is a small number.
 Let $\xi(x)$ be a lift of a ${\cal C}^{\infty}$ 1-periodic function satisfying
 \begin{equation}\label{phi0'}
\xi(x)=\left\{\begin{array}{ll} \xi_{01}(x),& \ \ |x|\le \delta;\\
 -\xi_{02}(x)\ ({\rm\ or\ } \pi+\xi_{02}(x)),& \ \ |x-\frac12|\le
 \delta. \end{array} \right.
\end{equation}
\noindent (b)  $\forall\ |x({\rm mod}\ 1)|>\delta\ \ {\rm and}\  |(x-1/2)({\rm mod}\ 1)|>\delta,$
\ $|\xi(x)({\rm mod}\ \pi)|> e^{-\frac{1}{\delta^a}}$.
\\
Define $\xi_N(x)=\xi(x)$ and $B_N(x)=\Lambda \cdot R_{\frac\pi 2-\xi_N(x)}$.

We restate Lemma 5.1 in \cite{WY} as follows:
\begin{Lemma}\label{L5.1}
For each $n\ge N$, there exist a $g_n(x)\in{\cal C}^{\infty}$ be a $1$-periodic function such that
  $$
  g_n(x):\ \ \left\{ \begin{array}{ll}=1,&\quad x\in\frac{I_n}{10},\\
  \\
  \in [0, 1],&\quad x\in I_n\backslash\frac{I_n}{10}\\
  \\
 =0,&\quad x\in \SS^1\backslash {I_n}
 \end{array}
 \right.
 $$
and
 \beq\label{dfdx}
 \left|\frac{d^rg_n(x)}{dx^r}\right|\le q_n^{3r},\quad 0\le r\le
[q_n^{\frac{1}{10}}].
 \eeq
\end{Lemma}
\vskip 0.3cm
Using the same argument as that in  finite smooth case, we have that
for any $x\in I_N$, $\|B^{r_N^+(x)}_N(x)\|\ge \l_N^{r_N^+(x)}$ and
\beq\label{difference0'}
 |\phi_{B_N,r_N}(x)+\psi_{B_N,-r_N}(x)-\frac{\pi}{ 2}-\xi_0(x)|_{{\cal C}^1}\le \l_N^{-1}
 \eeq
for $x\in I_N$.

 Define a $1$-periodic function $e_N(x)\in{\cal C}^{\infty}$ such that $e_N(x)=-(\phi_{B_N,r_N}(x)+\psi_{B_N,-r_N}(x)-\frac{\pi}{ 2}-\xi_0(x))$ for $x\in I_N$.
\vskip 0.3cm

 Let  $\hat{e}_N(x)=e_N(x)\cdot g_N(x)$ and  $\xi_{N+1}(x)=\xi_N(x)+\hat{e}_N(x)$ for $x \in \SS^1$.
Define $B_{N+1}(x)=\Lambda \cdot R_{\frac\pi2-\xi_{N+1}(x)}$. Obviously,
$B_{N+1}(x)=B_N(x)\cdot R_{-\hat{e}_N(x)}$. Then for any $x\in I_N$, $\|B^{r_N^+(x)}_{N+1}(x)\|\ge \l_N^{r_N^+(x)}$ and
$\phi_{B_{N+1},r_N}(x)+\psi_{B_{N+1},-r_N}(x)=\phi_{B_N,r_N}(x)+\psi_{B_N,-r_N}(x)-\hat{e}_N(x)$, which implies
$\phi_{B_{N+1},r_N}(x)+\psi_{B_{N+1},-r_N}(x)-\frac\pi2=\xi_0(x)$ on $\frac{I_N}{10}$. (\ref{difference0'})
implies that $|\hat{e}_N(x)|_{{\cal C}^1}\le \l_N^{-1}$ in
$I_N$. Thus we have $|\phi_{B_{N+1},r_N}(x)+\psi_{B_{N+1},-r_N}(x)-\frac\pi2|\ge \frac{1}{2}\cdot
 e^{-(10\cdot q_N^2)^{a}}$ on $I_N\backslash \frac{I_N}{10}$.\\



For any $n\ge N$, define a $1$-periodic function $e_{n}(x)\in{\cal C}^{\infty}$
such that
 \[e_{n}(x)=(\phi_{B_n,r_n}(x)+\psi_{B_n,-r_n}(x))-(\phi_{B_n,r_{n+1}}(x)+\psi_{B_n,-r_{n+1}}(x))
 \quad  x\in I_n.\]

  Define $\hat{e}_{n}(x)=e_{n}(x)\cdot g_{n}(x)$,  $\xi_{n}(x)=
  \xi_{n-1}(x)+\hat{e}_{n}(x)$ and $B_{n}(x)=\Lambda
\cdot R_{\frac\pi2-\xi_{n}(x)}$. Obviously, $B_n(x)=B_{n-1}(x)\cdot
R_{-\hat{e}_n(x)}$. Then we obtain (\ref{phin}), (\ref{difference-potential}) of Proposition \ref{WY2} and
  $$\begin{array}{ll}
   |\phi_{B_n,r_n}(x)+\psi_{B_n,-r_n}(x)-\frac{\pi}{ 2}|&=e^{-|x|^{-a}}({\rm\ or\ } e^{-|x-1/2|^{-a}}), \quad x\in \frac{I_{n,i}}{10}, i=1, 2\\
   |\phi_{B_n,r_n}(x)+\psi_{B_n,-r_n}(x)-\frac{\pi}{ 2}|&\ge \frac{1}{2}\cdot e^{-(10\cdot q_{n}^2)^{a}},\quad x\in I_n\backslash \frac{I_n}{10}.
  \end{array}
  $$
 \vskip 0.3cm
 \noindent From (\ref{phin}), one easily sees that $B_N(x),B_{N+1}(x), \cdots, $
 is ${\cal C}^1$-convergent to some $D_{\infty}(x)$.
 Furthermore, from (\ref{difference-potential}),
 the Lyapunov exponent of $D_{\infty}(x)$ has a lower bound
 $\log\l_{\infty}>(1-\e)\log \l$.\\

 In the following, we will prove that $B_N(x), B_{N+1}(x), \cdots, $
 is also convergent to $D_{\infty}(x)$ in ${\cal C}^{\infty}$-topology.

\vskip 0.3cm
 \begin{Lemma}
$B_N(x), B_{N+1}(x), \cdots, $
 is also convergent to $D_{\infty}(x)$ in ${\cal C}^{\infty}$-topology.
 \end{Lemma}
 \Proof
 It is equivalent to prove that $\xi_n(x),\ n=N, N+1, \cdots$ is ${\cal C}^{\infty}$-convergent. From the definition of $\xi_n(x)$, we have
 $\xi_{n}(x)-\xi_{n-1}(x)=\hat{e}_n(x)$. From the definition of $\hat{e}_n(x)$, it is sufficient to estimate
 $e_n(x)$ and $g_n(x)$. Since ${e}_n(x)$ is determined by
 $\phi_{B_n,r_n}(x)-\phi_{B_n,r_{n+1}}(x)$ and $\psi_{B_n,-r_n}(x)-\psi_{B_n,-r_{n+1}}(x)$,  with the help of Lemma \ref{iterative-estimate-derivative}, we have
 $$
\left|\frac{d^re_n(x)}{dx^r}\right|\le C(r)\cdot\l_n^{-q_{n-1}},\quad 0\le
r\le [q_{n-1}^{\frac{1}{10}}].
$$
 Note that $C(r)$ is independent of $n$. Thus for any fixed $R\in {\mathbb N}$, we can choose $n$ large enough such that $C(r)\le \l_n^{\frac12q_{n-1}}$ for any $r\le  R$. This together with (\ref{dfdx}) ends the proof.
 \hfill{}  \qedbox
\vskip 0.4cm

\noindent {\it Construction of $C_k(x)$}\quad
Next we will construct the sequence $C_k(x) (k=N, N+1, \cdots)$ with $L(C_k)=0$ such that  ${\cal C}^{\infty}$ converge to $D_{\infty}$.

 Consider the sub-interval $[0, \frac{1}{{q_{n_{i-1}}^2}}]$ of $I_{n_{i-1}}$. Define $0<c
<\tilde{c}< d<\frac{1}{{q_{n_{i-1}}^2}}$ such that
$|c|=(2\delta_0\cdot {r_{n_{i-1}}\cdot\log\underline\mu}_{n_{i-1}})^{-1/a}, |\tilde{c}|=M^2\cdot |c|, {d}=\frac{1}{2{q_{n_{i-1}}^2}}$.
Let $n_i$ be sufficiently large such that $I_{n_i}\nsupseteq [0, c]$.

 Define
$$
\bar{e}_i(x)=\left\{\begin{array}{ll}
e^{-|c|^{-a}}-(\phi_{ {A}_{i-1},r_{n_{i-1}}}(x)+\psi_{{A}_{i-1},-r_{n_{i-1}}}(x)),& x\in [\tilde{c}, {d}],\\
0, & x\notin [c, \frac{1}{q_{n_{i-1}}^2}],\\
\bar{h}_i^{\pm}(x) & x\in [c, \tilde{c}]\cup [d, \frac{1}{{q_{n_{i-1}}^2}}],
\end{array}
\right.
$$
where $\bar{h}_i^{\pm}(x)$ is of a ${\cal C}^{\infty}$ connection   between the parts in $[0, c]$ and $[\tilde{c}, {d}]$ as well as between the part in $[\tilde{c}, {d}]$ and the end point $\frac{1}{{q_{n_{i-1}}^2}}$ of $I_{n_i}$. Then similar to Lemma
\ref{L5.1}, we have  $$
 \left|\frac{d^r\bar{h}_i(x)}{dx^r}\right|\le C(r)\cdot q_{n_i}^{3r},\quad 0\le r\le
[q_{n_i}^{\frac{1}{10}}].
$$
Thus the ${\cal C}^{\infty}$-convergence of $C_k(x)$ is similar to the above argument. The remain part of the proof is same as Section \ref{s4}.

\appendix

\section {Product of hyperbolic matrices}

\vskip 0.8cm
{\rm  Let $A$ be a hyperbolic $SL(2, R)$-matrix, i.e., $\|A\|>1$.
It is know that   $A$ can be written uniquely as
$A=R_{\psi}\cdot \L_A\cdot R_{\phi}$ with $\L_A=diag(\|A\|, \|A\|^{-1})$. It is known that $-\phi$ is the most expanded direction of $A$ and
$\psi$ is the most contracted direction of $A^{-1}$.

 For two hyperbolic matrices $A=R_{\psi_A}\cdot \L_A\cdot R_{\phi_A}, B=R_{\psi_B}\cdot \L_B\cdot R_{\phi_B}$  with big norms, let $BA=R_{\psi_{BA}}\cdot \L_{BA}\cdot R_{\phi_{BA}}$. We firstly investigate how $\phi_{BA}, \psi_{BA}$ and $\|BA\|$ depend on $A$ and $B$.

\begin{Lemma}\label{L2.1}
Let $A, B$ be hyperbolic $SL(2, \mathbb{R})$ cocycles and $\theta=\phi_B+\psi_A$. Then it holds that
$\frac{1}{4}N(\|A\|, \|B\|, \t)\le \|BA\|^2\le N(\|A\|, \|B\|, \t)$, where $N(\|A\|, \|B\|, \t)=(\|A\|^2\|B\|^2+\|A\|^{-2}\|B\|^{-2})\cdot \cos^2\theta+(\|A\|^2\|B\|^{-2}+\|A\|^{-2}\|B\|^2)\cdot\sin^2\theta$.
\end{Lemma}
\begin{Proof}
For any $SL(2, \mathbb{R})$ matrix $A=(a_{ij})_{2\times 2}$, it is known that $\frac14 \sum_{i,j}a_{ij}^2\le \|A\|^2\le \sum_{i,j}a_{ij}^2$.

It is easy to see that
$$\begin{array}{ll}\|BA\|&=\left\|\left(\begin{array}{ll}\|B\| & 0\\ 0& \|B\|^{-1}\end{array}\right)\cdot
\left(\begin{array}{ll}\cos \theta & -\sin\theta\\ \sin\theta& \cos\theta\end{array}\right)
\cdot \left(\begin{array}{ll}\|A\| & 0\\ 0& \|A\|^{-1}\end{array}\right)\right\|\\
\\
&=\left\|\left(\begin{array}{ll}\|A\|\|B\|\cos\theta & -\|A\|^{-1}\|B\|\sin\theta\\ \|A\|\|B\|^{-1}\sin\theta& \|A\|^{-1}\|B\|^{-1}\cos\theta\end{array}\right)\right\|.\end{array}$$
It thus implies the conclusion. \quad \qedbox
\end{Proof}

\begin{Lemma}\label{L2.2}Let  $\phi=\phi_A-\phi_{BA}$, $\psi=\psi_{BA}-\psi_B$. Assume $\t\in [0,\pi)$. Then \begin{equation}\label{phi00}
\phi(\|A\|, \|B\|, \t)=\left\{\begin{array}{ll}0, &\quad{\rm for} \ \ \t=0\\
\\
-\frac 12\left(\frac\pi2- \tan^{-1}(a\cot \t+b\tan \t)\right), &\quad{\rm for} \ \ 0<\t< \frac\pi2\\
\\
\frac\pi2-\frac 12\left(\frac\pi2- \tan^{-1}(a\cot \t+b\tan \t)\right), &\quad{\rm for} \ \ \frac\pi2<\t<\pi {\rm\ if}\ b\ge 0\\
\\
-\frac\pi2-\frac 12\left(\frac\pi2- \tan^{-1}(a\cot \t+b\tan \t)\right), &\quad{\rm for} \ \ \frac\pi2<\t<\pi {\rm\ if}\ b< 0\\
\\
0, &\quad{\rm for} \ \ \t=\frac\pi2 {\rm\ if}\ b\ge 0\\
\\
-\frac\pi2, &\quad{\rm for} \ \ \t=\frac\pi2 {\rm\ if}\ b<0,
 \end{array} \right.
\end{equation}
where
$$
 a=\frac{\|A\|^{2}\|B\|^2-\|A\|^{-2}\|B\|^{-2}}{2(\|B\|^2-\|B\|^{-2})}, \quad
b=\frac{\|A\|^{2}\|B\|^{-2}-\|A\|^{-2}\|B\|^{2}}{2(\|B\|^2-\|B\|^{-2})}.
$$
 Similarly,
\begin{equation}\label{psi0}
\psi(\|A\|, \|B\|, \t))=\left\{\begin{array}{ll}0, &\quad{\rm for} \ \ \t=0\\
\\
-\frac 12\left(\frac\pi2- \tan^{-1}(a'\cot \t-b'\tan \t)\right), &\quad{\rm for} \ \ 0<\t< \frac\pi2\\
\\
\frac\pi2-\frac 12\left(\frac\pi2- \tan^{-1}(a'\cot \t-b'\tan \t)\right), &\quad{\rm for} \ \ \frac\pi2<\t<\pi {\rm\ if}\ b'\ge 0\\
\\
-\frac\pi2-\frac 12\left(\frac\pi2- \tan^{-1}(a'\cot \t-b'\tan \t)\right), &\quad{\rm for} \ \ \frac\pi2<\t<\pi {\rm\ if}\ b'< 0\\
\\
0, &\quad{\rm for} \ \ \t=\frac\pi2 {\rm\ if}\ b'\ge 0\\
\\
-\frac\pi2, &\quad{\rm for} \ \ \t=\frac\pi2 {\rm\ if}\ b'< 0,
 \end{array} \right.
\end{equation}
where
$$
 a'=\frac{\|A\|^{2}\|B\|^2-\|A\|^{-2}\|B\|^{-2}}{2(\|A\|^2-\|A\|^{-2})}, \quad
b'=\frac{\|A\|^{2}\|B\|^{-2}-\|A\|^{-2}\|B\|^{2}}{2(\|A\|^2-\|A\|^{-2})}.
$$
\end{Lemma}
\begin{Proof}
Let
$$
\begin{array}{ll}
V(s)&=\left(\begin{array}{ll}\|B\| & 0\\ 0& \|B\|^{-1}\end{array}\right)\cdot
\left(\begin{array}{ll}\cos \theta & -\sin\theta\\ \sin\theta& \cos\theta\end{array}\right)
\cdot \left(\begin{array}{ll}\|A\| & 0\\ 0& \|A\|^{-1}\end{array}\right)\cdot \left(\!\!\!\!\begin{array}{lll}&\cos s& \\ &\sin s& \end{array}\!\!\!\!\right)\\
\\
&=\left(\begin{array}{ll}\|B\| & 0\\ 0& \|B\|^{-1}\end{array}\right)\cdot
\left(\begin{array}{ll}\cos\theta\cdot \|A\|\cdot \cos s-\sin\theta\cdot \|A\|^{-1}\cdot \sin s & \\ \sin\theta \|A\|\cdot\cos s+\cos\theta\cdot \|A\|^{-1}\cdot \sin s \end{array}\right)\\
\\
&=\left(\begin{array}{ll}\cos\theta\cdot \|A\|\|B\|\cdot \cos s-\sin\theta\cdot \|A\|^{-1}\cdot \|B\| \sin s & \\ \sin\theta \|A\|\|B\|^{-1}\cdot\cos s+\cos\theta\cdot \|A\|^{-1}\|B\|^{-1}\cdot \sin s \end{array}\right).
\end{array}$$
Thus
$$
\begin{array}{ll}
 |V(s)|^2
=&(\cos\theta \|A\|\|B\|)^2+(\sin^2\theta \|A\|^{-2}\|B\|^2-\cos^2\theta \|A\|^2\|B\|^2)\sin^2 s+\sin^2\theta
\|A\|^2\|B\|^{-2}\\
\\
&+(\cos^2\theta \|A\|^{-2}\|B\|^{-2}
-\sin^2\theta \|A\|^2\|B\|^{-2})\sin^2s
+2(\|B\|^{-2}-\|B\|^2)\sin\theta \cos\theta\sin s \cos s.
\end{array}
$$
Obviously $\frac{d}{ds}(|V(s)|^2)=0$ at $\phi$ since $|V(s)|^2$  attains its extreme at $\phi$,
a simple computation leads to
$$\begin{array}{ll}
& \left( ( \|A\|^{2}\|B\|^2-\|A\|^{-2}\|B\|^{-2})\cos^2 \t +(\|A\|^{2}\|B\|^{-2}-\|A\|^{-2}\|B\|^{2})\sin^2\t\right)\sin2\phi\\
&=-2(\|B\|^{2}-\|B\|^{-2})\sin2\theta\cos 2\phi.
\end{array}$$
Thus
$$
-\cot 2\phi=  \frac{\|A\|^{2}\|B\|^2-\|A\|^{-2}\|B\|^{-2}}{2(\|B\|^2-\|B\|^{-2})}\cot \t+\frac{\|A\|^{2}\|B\|^{-2}-\|A\|^{-2}\|B\|^{2}}{2(\|B\|^2-\|B\|^{-2})}\tan \t.
$$
With the help of the inequality $\frac{d^2}{ds^2}(|V(s)|^2)\le 0$, we obtain the unique $\phi$ corresponding the maximum $\|BA\| ^2$ of $|V(s)|^2$, which satisfies (\ref{phi00}).

(\ref{psi0}) is proved similarly.\qquad \qedbox
\end{Proof}

\vskip 0.4cm

Later we will see that both $\|A\|$ and $\|B\|$ are very big. Thus
$$\begin{array}{ll}
 a&=\frac{\|A\|^{2}\|B\|^2-\|A\|^{-2}\|B\|^{-2}}{\|B\|^2-\|B\|^{-2}}\sim \|A\|^{2},\quad\\
 \\
b&=\frac{\|A\|^{2}\|B\|^{-2}-\|A\|^{-2}\|B\|^{2}}{\|B\|^2-\|B\|^{-2}}\lesssim\max\{\|A\|^{-2},
\frac{\|A\|^{2}}{\|B\|^4}\}.\end{array}
$$

If $A, B$ are hyperbolic,
the functions $\phi(\|A\|, \|B\|, \theta), \psi(\|A\|, \|B\|, \theta)$ defined above  are continuous in all variables.  In the following, we estimate  the derivatives of $\phi$ and $\psi$  with respect to $\theta,\ \|A\|$ and $\|B\|$.
\begin{Lemma}\label{pa phi pa t}
It holds that
\beq\label{paphipat0}
|\phi (mod\ \pi)|\le C(0)\cdot \|A\|^{-2}\cdot |\t-\frac{\pi}{2}|
^{-1}
\eeq
and
\beq\label{papsipat0}
|\psi (mod\ \pi)|\le C(0)\cdot \|B\|^{-2}\cdot |\t-\frac{\pi}{2}|
^{-1}.
\eeq
Suppose $|\t-\frac\pi2|^{-1}\ll \|A\|^2$. Then, for $i\ge 1$, we have that \beq\label{paphipat}
|\frac{\pa^i\phi }{\pa\t^i}|\le C(i)\cdot \|A\|^{-2}\cdot |\t-\frac{\pi}{2}|
^{-i-1},
\eeq
\beq\label{1}\begin{array}{ll}
|\frac{\pa^i\phi}{\pa \|A\|^i}|&\le
 C(i)\cdot \|A\|^{-2}\cdot \|A\|^{-i}\cdot |\t-\frac{\pi}{2}|^{-1},
\end{array}
\eeq
and
\beq\label{2}
|\frac{\pa^i\phi}{\pa \|B\|^i}|\le C(i)\cdot \|A\|^{-2}\cdot \|B\|^{-i}\cdot |\t-\frac{\pi}{2}|^{-1}.
\eeq
More generally, for $i+j+k\ge 1$, we have
\beq\label{ij}
|\frac{\pa^{i+j+k}\phi}{\pa \t^i\pa\|A\|^j\pa\|B\|^k}|\le C(i,j,k)\cdot |\t-\frac\pi2|^{-i-1}\|A\|^{-2-j}\cdot \|B\|^{-k};
\eeq

Similarly, suppose $|\t-\frac\pi2|^{-1}\ll \|B\|^2$. Then we have
\beq\label{ijk}
|\frac{\pa^{i+j+k}\psi}{\pa \t^i\pa\|A\|^j\pa\|B\|^k}|\le C(i,j,k)\cdot |\t-\frac\pi2|^{-i-1}\|A\|^{-j}\cdot \|B\|^{-2-k}.
\eeq

\end{Lemma}
\begin{Proof} To prove (\ref{paphipat0}), it is sufficient to consider the situation $\t\approx \frac\pi2$. We only consider the case $0\le \t\le \frac\pi2$ since the proof for the other cases is similar. From the fact
$\lim_{x\rightarrow \infty}\frac{\frac\pi2-\arctan x}{x^{-1}}=1$ and the definition of $a$, we have
$|\phi|\le C(0) \cdot a^{-1}\cdot |\t-\frac\pi2|^{-1}\le C(0) \cdot \|A\|^{-2}\cdot |\t-\frac\pi2|^{-1}$. Thus we obtain
(\ref{paphipat0}). We can obtain (\ref{papsipat0}) similarly.

for $i\ge1$, from the definition of $\phi$, we have
$$
\frac{\pa^i \phi }{\pa\t^i}=-\frac12\sum_{l_1+\cdots +l_k=i}\frac{d^{k-1}(\frac{1}{1+f^2})}{df^{k-1}}\cdot\frac{\pa^{l_1}f}
{\pa\t^{l_1}}\cdots \frac{\pa^{l_k}f}{\pa\t^{l_k}},$$
where  $f(\|A\|, \|B\|, \theta)=a\cot \t+b\tan \t$.
To estimate $\frac{\pa^i \phi}{\pa\t^i}$, we have that
$$
\begin{array}{ll}
|\frac{\pa^{l_s}f}{\pa\t^{l_s}}|&=|\frac{\pa^{l_s}}{\pa\t^{l_s}}(a\cot \t+b \tan\t)|\le |a|
\cdot |\cot^{(l_s)}(\t)|+|b|\cdot |\tan^{(l_s)}(\t)|.
\end{array}
$$
By a direct computation, we have
$$
|\tan^{(l_s)}\t|=|(\cos^{-2}\t)^{(l_s-1)}|\le |\sum_{\kappa_1+\cdots +\kappa_t=l_s-1}
\cos^{-(2+t)}\t\cdot \cos^{(\kappa_1)}\t\cdots \cos^{(\kappa_t)}\t|
$$
and
$$
|\cot^{(l_s)}\t|=|(\sin^{-2}\t)^{(l_s-1)}|\le |\sum_{\kappa_1+\cdots +\kappa_t=l_s-1}
\sin^{-(2+t)}\t\cdot \sin^{(\kappa_1)}\t\cdots \sin^{(\kappa_t)}\t|.
$$

From the condition $|\t-\frac\pi2|^{-1}\ll \|A\|^2$ and the fact that the signs of $\|B\|^2\cot \t$ and $\|B\|^{-2}\tan \t$ are the same, we have
\beq\label{3}
|\frac{\pa^{l_s}f}{\pa\t^{l_s}}|\le  C(l_s)\cdot(|a|\cdot|\t-\frac{\pi}{2}|^{-(l_s-1)}+|b|\cdot |\t-\frac{\pi}{2}|^{-(l_s+1)})\le C(l_s)\cdot |f| \cdot |\t-\frac\pi2|^{-l_s}.
\eeq
On the other hand, we have
$$
\left|\frac{d^{k-1}(\frac{1}{1+f^2})}{df^{k-1}}\right|\lesssim
|f|^{-k-1}\quad {\rm if}\ k\ge 1.
$$
Thus from $|f|\gtrsim \|A\|^2\cdot |\cot \t|$ we obtain
$$
|\frac{\pa^i\phi}{\pa\t^i}|\le C(i)|\t-\frac{\pi}{2}|^{-i}\cdot\frac{1}{|f|}
\le C(i)\|A\|^{-2}|\t-\frac{\pi}{2}|
^{-i-1}.
$$

Next we estimate \beq\label{phiA}\begin{array}{ll}
|\frac{\pa^i\phi}{\pa \|A\|^i}|&\le \sum_{l_1+\cdots +l_k=i}  |\frac{d^{k-1}(\frac{1}{1+f^2})}{df^{k-1}}|  \cdot
|\frac{\pa^{l_1}f}{\pa \|A\|^{l_1}}|\cdots |\frac{\pa^{l_k}f}{\pa \|A\|^{l_k}}| \quad l_j\ge 1,\ 1\le j\le k\\
\\
&\le \sum_{l_1+\cdots +l_k=i} |f|^{-k-1} \cdot
|\frac{\pa^{l_1}f}{\pa \|A\|^{l_1}}|\cdots |\frac{\pa^{l_k}f}{\pa \|A\|^{l_k}}|.
\end{array}\eeq

It is easy to see that $|f|\sim |a|\cdot |\cot \t|$ with the condition $|\t-\frac\pi2|^{-1}\ll \|A\|^2$. We also have
$$
|\frac{\pa^{l_s}f}{\pa \|A\|^{l_s}}|\le |\cot \t||\frac{\pa^{l_s}a}{\pa \|A\|^{l_s}}|
+|\tan \t||\frac{\pa^{l_s}b}{\pa \|A\|^{l_s}}|.
$$
By a direct computation, we obtain
$$
|\frac{\pa^{l_s}a}{\pa \|A\|^{l_s}}|=\left|\frac{\pa^{l_s}}{\pa \|A\|^{l_s}}\left(
\frac{\|A\|^{2}\|B\|^2-\|A\|^{-2}\|B\|^{-2}}{\|B\|^2-\|B\|^{-2}}\right)
\right|\le C(l_s)\cdot |a|\cdot \|A\|^{-l_s}
$$
and
$$\begin{array}{ll}
|\frac{\pa^{l_s}b}{\pa \|A\|^{l_s}}|=\left|\frac{\pa^{l_s}}{\pa \|A\|^{l_s}}\left(
\frac{\|A\|^{2}\|B\|^{-2}-\|A\|^{-2}\|B\|^{2}}{\|B\|^2-\|B\|^{-2}}\right)
\right|\le C(l_s)\cdot(\|A\|^{-2}+\frac{\|A\|^{2}}{\|B\|^4})\cdot \|A\|^{-l_s}.
\end{array}
$$
Thus we have
\beq\label{fA}\begin{array}{ll}
|\frac{\pa^{l_s}f}{\pa \|A\|^{l_s}}|&\le C(l_s)\cdot \left\{ |\t-\frac{\pi}{2}|\|A\|^{2-l_s}
+|\t-\frac{\pi}{2}|^{-1}\cdot \|A\|^{-l_s}\cdot (\|A\|^{-2}+\frac{\|A\|^{2}}{\|B\|^4})\right\}\\
\\
&\le C(l_s)\cdot |f|\cdot \|A\|^{-l_s}.
\end{array}
\eeq
With the fact that $|f|\gtrsim \|A\|^2\cdot |\cot \t|$, it follows that
\beq\label{f-2A}\begin{array}{ll}
|f|^{-2}\cdot|\frac{\pa^{l_s}f}{\pa \|A\|^{l_s}}|&\le
{C(l_s)} |f|^{-1}\cdot \|A\|^{-l_s}
\\
\\
&\le C(l_s)\cdot \|A\|^{-2-l_s}\cdot |\t-\frac{\pi}{2}|^{-1}.
\end{array}
\eeq

Combining (\ref{phiA}), (\ref{fA}) with (\ref{f-2A}), we obtain (\ref{1}).
Similarly, we have (\ref{2}) and (\ref{ij}).
The estimates for $\psi$ can be proved similarly.
\end{Proof}\qquad \qedbox



\section{Proof of Lemma \ref{iterative-estimate-derivative}}
In this section, we first give estimates on most contracted and expanded directions of the product of hyperbolic blocks.
Then we will give the proof of Lemma \ref{iterative-estimate-derivative}.

Let $A(x), B(x), \theta(x), \phi(x)$ and $\psi(x)$ be defined as in Lemma \ref{L2.1} and \ref{L2.2}.


\vskip 0.4cm
\begin{Lemma}\label{l_{n_1+n_2}}
Let $|\t-\frac{\pi}{2}|^{-1}\ll \|A\|^2, \|B\|^2$. Suppose that, for any $i\ge 0$,
\beq\label{frac{dl_{n_1 or n_2}}{dx}}
\left|\frac{d^i\|A\|}{dx^i}\right|\le C(i)\cdot \|A\|\cdot |\t-\frac{\pi}{2}|^{-i-1},  \quad \left|\frac{d^i\|B\|}{dx^i}\right|\le C(i)\cdot \|B\|\cdot |\t-\frac{\pi}{2}|^{-i-1},\ \
|\frac{d^i\t}{dx^i}|\le C(i)\cdot|\t-\frac{\pi}{2}|^{-i-1}.
\eeq
Then we have
\beq\label{frac{dl_{n_1+n_2}}{dx}}
\begin{array}{ll}
&|\frac{d^i\phi}{dx^i}|\le C(i)\cdot |\t-\frac{\pi}{2}|^{-i-1}\cdot  \|A\|^{-2}, \\
\\
&|\frac{d^i\psi}{dx^i}|\le C(i)\cdot |\t-\frac{\pi}{2}|^{-i-1}\cdot  \|B\|^{-2},
\\
\\
&\left|\frac{d^i\|BA\|}{dx^i}\right|\le C(i)\cdot \|BA\|\cdot |\t-\frac{\pi}{2}|^{-i-1}.\end{array}
\eeq
\end{Lemma}
\begin{Proof}
 From Lemma \ref{pa phi pa t} and
 \beq\label{dphidx}\begin{array}{ll}&\left|\frac{d^i\phi}{dx^i}\right|=\sum_{t_1+\cdots+t_{i_1}+s_1+\cdots+s_{i_2}+j_1+\cdots +j_{i_3}=i}\left|\frac{\pa ^{i_1+i_2+i_3}\phi}{
\pa \|A\|^{i_1}\cdot \pa \|B\|^{i_2}\pa \t^{i_3}}\cdot \frac{d^{t_1}\|A\|}{d x^{t_1}}\cdots \frac{d^{t_{i_1}}\|A\|}{d x^{t_{i_1}}}\right|\\
\\
&\cdot  \left|\frac{d^{s_1}\|B\|}{d x^{s_1}}\cdots \frac{d^{s_{i_2}}\|B\|}{d x^{s_{i_2}}}
\cdot \left(\frac{d^{j_1}\t}{d x^{j_1}}\right)\cdots
\left(\frac{d{j_{i_3}}\t}{d x^{j_{i_3}}}\right)\right|
,\end{array}
\eeq
we have
\beq\label{8}\begin{array}{ll}\left|\frac{\pa ^{i_1+i_2+i_3}\phi}{
\pa \|A\|^{i_1}\cdot \pa \|B\|^{i_2}\cdot\pa \t^{i_3}}\right|
\le C(i_1,i_2,i_3)\cdot |\t-\frac{\pi}{2}|^{-(1+i_3)}\cdot \|A\|^{-i_1-2}\cdot
\|B\|^{-i_2}.
\end{array}
\eeq
Then from (\ref{frac{dl_{n_1 or n_2}}{dx}}), (\ref{dphidx}), (\ref{8}), we have
$$
\left|\frac{d^i\phi}{dx^i}\right|\le C(i)\cdot \|A\|^{-2}\cdot |\t-\frac{\pi}{2}|^{-(i+1)},
$$
thus the first inequality of (\ref{frac{dl_{n_1+n_2}}{dx}}) is proved.
In  the above inequality, we used the fact
$$
\left|\frac{\pa f}{\pa \|A\|}\cdot \frac{\pa^j \|A\|}{\pa x^j}\right|\le |\t-
\frac{\pi}{2}|^{-(j+1)}\cdot |f|.
$$

The second inequality is proved similarly. Now we prove the third inequality.
By a direct computation, we have
\beq\label{frac{pa l}{pa phi}}\frac{\pa^i \|BA\|}{\pa \phi^i}=\frac{\pa^i(g^{\frac 1 2})}{\pa\phi^i}
=\sum_{l_1+\cdots +l_k=i}(g^{\frac12})^{(k)}\cdot \frac{\pa^{l_1}g}{\pa \phi^{l_1}}
\cdots \frac{\pa^{l_k}g}{\pa \phi^{l_k}},
\eeq
where
\beq\label{gg}\begin{array}{ll}
&g=g^2_1+g^2_2, \quad g_1=\|A\|\cdot \|B\|\cdot \cos\t \cos \phi-
\|A\|^{-1}\cdot \|B\|^{-1}\sin\t\sin\phi,\\
\\
& g_2=\|A\|\cdot
\|B\|^{-1}\cdot \sin\t\cdot \cos \phi+\|A\|^{-1}\cdot \|B\|^{-1}
\cdot \cos\t\cdot \sin\phi.\end{array}
\eeq
It is not difficult to see that $|(g^{\frac12})^{(k)}|\le C(k)\cdot g^{\frac12-k}$.
From the definition of $g$, we have
$$
\left|\frac{\pa^{l_s}g}{\pa\phi^{l_s}}\right|\le \left|\frac{\pa^{l_s}(g_1^2)}{\pa\phi^{l_s}}\right|+
\left|\frac{\pa^{l_s}(g_2^2)}{\pa\phi^{l_s}}\right|.$$
From
$$
\begin{array}{ll}
\left|\frac{\pa^{l_s}(g_1^2)}{\pa\phi^{l_s}}\right|\le \sum_{l_{s,1}+l_{s,2}=l_s}\left|\frac{\pa^{l_{s,1}} g_1}{\pa \phi^{l_{s,1}}}\right|\cdot \left|\frac{\pa^{l_{s,2}} g_1}{\pa \phi^{l_{s,2}}}\right|,\end{array}
$$
it follows that
$$\begin{array}{ll}
\left|\frac{\pa^{l_{s,1}} g_1}{\pa \phi^{l_{s,1}}}\right|&\le \|A\|\|B\||\cos\t|\cdot
|\cos(\phi+\frac{\pi}{2}\cdot l_{s,1})|+\|A\|^{-1}\|B\||\sin \t|\cdot |\sin(\phi+\frac{\pi}{2}\cdot l_{s,1})|\\
\\
&\le \|A\|\cdot \|B\|\cdot |\cos\t|\lesssim\|BA\|.
\end{array}
$$
Then we obtain
\beq \label{g^2_1}
\left|\frac{\pa^{l_s}(g_1^2)}{\pa\phi^{l_s}}\right|\lesssim \|BA\|^2.\eeq
Similarly, we have
\beq \label{g^2_2}
\left|\frac{\pa^{l_s}(g_2^2)}{\pa\phi^{l_s}}\right|\lesssim \|BA\|^2.\eeq
Combining (\ref{frac{pa l}{pa phi}}) with (\ref{g^2_1}), (\ref{g^2_2}), we then have
\beq\label{4}
\left|\frac{\pa^{i}\|BA\|}{\pa\phi^{i}}\right|\le C(i)\cdot \max_{k\le i}(\|BA\|^{2k}\cdot g^{\frac12 -k})
=C(i)\cdot \|BA\|.
\eeq
Similarly, it holds that
\beq \label{5}\left|\frac{\pa^i\|BA\|}{\pa \t^i} \right|\le C(i)\cdot  \max_{k\le i} \left(g^{\frac12-k}\cdot (\|A\|\|B\|)^{2k}|\cos\t|^k\right)\le C(i)\cdot |\t-\frac{\pi}{2}|^{1-i}\eeq
and
\beq \label{67}\left|\frac{\pa^i\|BA\|}{\pa \|A\|^i} \right|\le C(i)\cdot \|BA\|\cdot \|A\|^{-i},\quad
\left|\frac{\pa^i\|BA\|}{\pa \|B\|^i} \right|\le C(i)\cdot \|BA\|\cdot \|B\|^{-i}.\eeq

Similar to (\ref{4})-(\ref{67}),
 we have
$$\begin{array}{ll}
\left|\frac{\pa^{i_1+\cdots +i_4}\|BA\|}{\pa \|A\|^{i_1}\cdot
\pa \|B\|^{i_2}\cdot\pa^{i_3}\phi\cdot\pa^{i_4}\t}\right|
&\le C(i)\cdot \|BA\|\cdot \|A\|^{-i_1}\cdot \|B\|^{-i_2}\cdot |\t-\frac{\pi}{2}|^{-i_4}.
\end{array}
$$
Combining with (\ref{l_{n_1+n_2}}), the first inequality in (\ref{frac{dl_{n_1+n_2}}{dx}}) and the fact
\[
\frac{d^i\|BA\|}{dx^i}=\sum\frac{\pa^{i_1+\cdots +i_4}\|BA\|}{\pa \|A\|^{i_1}\cdot
\pa \|B\|^{i_2}\cdot\pa^{i_3}\phi\cdot\pa^{i_4}\t}\cdot
\frac{\pa^{j_{1,1}}\|A\|}{\pa x^{j_{1,1}}}\cdots \frac{\pa^{j_{i_1,1}}\|A\|}{\pa x^{j_{i_1,1}}}
\cdots \frac{\pa^{j_{1,4}}\t}{\pa x^{j_{1,4}}} \cdots \frac{\pa^{j_{i_4,4}}\t}{\pa x^{j_{i_4,4}}},\]
 we prove the third inequality in (\ref {frac{dl_{n_1+n_2}}{dx}}).\qquad \qedbox
\end{Proof}
\vskip 0.4cm

\noindent
{\it Proof of Lemma  \ref{iterative-estimate-derivative}}.
For any $x\in I_{k+1}$, let $r_k(x):=r_{k,0}(x)<r_{k,1}(x)<\cdots <r_{k,s(x)}(x):=r_{k+1}(x)$ such that
$T^{r_{k,j}(x)}x\in I_k,\ 0\le j\le s(x)\le C(M)$. Consider $A^{r_{k,0}(x)+r_{k,1}(x)}(x)=A^{r_{k,1}(x)}(T^{r_{k,0}(x)}(x))\cdot
A^{r_{k,0}(x)}(x)$.

Let $A^{r_{k,0}(x)}(x):=R_{\psi_k^-(x)}\cdot \Lambda_{k^-}(x)\cdot R_{\phi_k^-(x)}$ and $A^{r_{k,1}(x)}(T^{r_{k,0}(x)}x):=R_{\psi_k^+(x)}\cdot \Lambda_{k^+}(x)\cdot R_{\phi_k^+}(x)$.
Then $$A^{r_{k,0}(x)+r_{k,1}(x)}(x)=R_{\psi_k^+(x)}\cdot \Lambda_{k^+}(x)\cdot R_{\phi_k^++\psi_k^-}(x)\cdot \Lambda_{k^-}(x)\cdot R_{\phi_k^-(x)}:=R_{\psi_{k+1,1}}(x)\cdot \Lambda_{k+1,1}(x)\cdot R_{\phi_{k+1,1}}(x).$$
Since $T^{r_{k,j}(x)}x\in I_k\backslash I_{k+1}$ for $j<s(x)$, it holds that $|\phi_k^++\psi_k^--\frac{\pi}{2}|\ge d_{k+1}$. From (\ref{norm-n}) and Lemma \ref{l_{n_1+n_2}}, we have
\beq\label{iterate-diff}
\left|\frac{d^i(\phi_{k+1,1}-\phi_k^-)}{dx^i} \right|\le C(i)\cdot d_k^{-i}\cdot \|A^{r_{k,1}(x)}(x)\|^{-2}.
\eeq
Similarly, it holds that
\beq\label{iterate-diff1}
\left|\frac{d^i(\psi_{k+1,1}-\psi_k^+)}{dx^i} \right|\le C(i)\cdot d_k^{-i}\cdot \|A^{r_{k,0}(x)}(x)\|^{-2}.
\eeq

In the above, we regard $\phi_{k+1,1}-\phi_k^-$  and $\phi_k^++\psi_k^-$ as $\phi$ and $\theta$ in Lemma \ref{l_{n_1+n_2}}, respectively.
Moreover, for each $x\in I_k$, $|\theta(x)-\frac\pi2|=|\phi_k^+(x)+\psi_k^-(x)-\frac{\pi}{2}|\ge d_{k+1}$ from the definition of $d_{k+1}$.
It implies that
$$\begin{array}{ll}
&\left|\frac{d^i\phi_{k+1,1}}{dx^i}\right|,\quad \left|\frac{d^i\psi_{k+1,1}}{dx^i}\right|\le C(i)\cdot d_k^{-i}+\left|
\frac{d^i\phi_k^-}{dx^i}\right|+\left|\frac{d^i\psi_k^+}{dx^i}\right|\\
\\
&\le C(i)\cdot (d_k^{-i}+d_{k-1}^{-i})\le C(i)\cdot d_k^{-i}.
\end{array}
$$
The last inequality is obtained from ${\rm (1)}_k$.

From Lemma \ref{l_{n_1+n_2}} we obtain $\left|\frac{d^i\|A^{r_{k,0}(x)+r_{k,1}(x)}(x)\|}{dx^i}\right|\le C(i)\cdot \|A^{r_{k,0}(x)+r_{k,1}(x)}(x)\|\cdot d_k^{-i}.$
Since $s(x)\le C(M)$, it follows from no more than $C(M)$-applications of the above argument that ${\rm (1)}_{k+1}$ and ${\rm (2)}_{k+1}$ hold true.


 Iterating (\ref{iterate-diff}) and (\ref{iterate-diff1}) less than $C(M)$ times, we get
$$
\left|\frac{d^i}{dx^i}(\phi_{A,r^+_{k+1}}(x)-\phi_{A,r^+_{k}}(x))\right|\le C(i)\cdot |\theta_k-\frac{\pi}{2}|^{-i}\cdot \|A^{r^+_k}\|^{-2} $$
and
$$
\left|\frac{d^i}{dx^i}(\psi_{A,-r^-_{k+1}}(x)-\psi_{A,-r^-_{k}}(x))\right|\le C(i)\cdot |\theta_k-\frac{\pi}{2}|^{-i}\cdot \|A^{r^-_k}\|^{-2}.
$$
 Since
$|\theta_k-\frac{\pi}{2}|\ge d_k$, we obtain (\ref{difference-derivative}).


\hfill{} \qedbox
\begin{Remark}
In the proof of Lemma  \ref{iterative-estimate-derivative}, it is not necessary that $s(x)$ is bounded by a constant. We make such an assumption only for the simplicity. And the condition that $\omega$ is of bound type is only used in  constructing $C_k(x)$.
\end{Remark}}

\bigskip
\footnotesize
\noindent\textit{Acknowledgments.}
The authors would like to thank S. Jitomirskaya for drawing their attention to the problem.
Y.Wang was  supported by NSFC (grant no.11271183), J. You was supported by NSFC (grant no. 11471155) and 973 projects of China (grant no. 2014CB340701).

\end{document}